\documentclass[twocolumn,amssymb,aps,secnumarabic,nofootinbib]{revtex4}
\usepackage{times,mathptmx,amsthm,amsmath}
\usepackage{pstricks,pst-node,graphics}
\newcommand{\half}{{\small 1/2}}
\newcommand{\odisk}{\pscircle[fillstyle=solid]}

\newcommand{\onex}{1.732}
\newcommand{\twox}{3.464}
\newcommand{\threex}{5.196}
\newcommand{\fourx}{6.928}

\theoremstyle{definition}

\psset{linewidth=.5pt}

\begin{document}
\title{Generalized Domino-Shuffling}
\author{James Propp ({\tt propp@math.wisc.edu})}
\thanks{Supported by grants from the National Science
Foundation and the National Security Agency.} 
\affiliation{Department of Mathematics \\
University of Wisconsin - Madison}
\date{September 30, 2001}
\begin{abstract}
\vspace{.1in}
The problem of counting tilings of a plane region
using specified tiles
can often be recast as 
the problem of counting (perfect) matchings 
of some subgraph of an Aztec diamond graph $A_n$,
or more generally 
calculating the sum of the weights of all the matchings,
where the weight of a matching is equal to the product
of the (pre-assigned) weights of the constituent edges
(assumed to be non-negative).
This article presents efficient algorithms 
that work in this context to solve three problems: 
finding the sum of the weights of the matchings
of a weighted Aztec diamond graph $A_n$; 
computing the probability that a randomly-chosen matching of $A_n$ 
will include a particular edge
(where the probability of a matching is proportional to its weight); 
and generating a matching of $A_n$ at random.
The first of these algorithms is equivalent to
a special case of Mihai Ciucu's cellular complementation algorithm \cite{Ci2}
and can be used to solve many of the same problems.
The second of the three algorithms is a generalization 
of not-yet-published work of Alexandru Ionescu, 
and can be employed to prove an identity 
governing a three-variable generating function
whose coefficients are all the edge-inclusion probabilities;
this formula has been used \cite{CEP}
as the basis for asymptotic formulas for these probabilities,
but a proof of the generating function identity
has not hitherto been published.
The third of the three algorithms 
is a generalization of the domino-shuffling algorithm
presented in \cite{EKLP}; 
it enables one to generate random ``diabolo-tilings of fortresses''
and thereby to make intriguing inferences 
about their asymptotic behavior.
\vspace{.1in}
\end{abstract}
\maketitle

\section{Introduction}

\subsection{Background}

Let $L$ be the rotated square lattice with vertex set 
$\{(i,j): i,j \in {\bf Z}, \ \mbox{$i+j$ is odd}\}$, 
and with an edge joining two vertices in the graph 
if and only if the Euclidean distance between the corresponding 
points in the plane is $\sqrt{2}$.
Define the {\it Aztec diamond graph\/} of order $n$
(denoted by $A_n$) 
as the induced subgraph of $L$ 
with vertices $(i,j$) satisfying $-n \leq i,j \leq n$, 
as shown below for $n=3$.
$$
\pspicture(-1.8,-1.8)(1.8,1.8)
\psset{xunit=.5cm,yunit=.5cm}
\pspolygon(-3,-2)(2,3)(3,2)(-2,-3)
\pspolygon(-3,0)(0,3)(3,0)(0,-3)
\pspolygon(-3,2)(-2,3)(3,-2)(2,-3)
\endpspicture
$$
A {\it partial matching\/} of a graph 
is a set of vertex-disjoint edges belonging to the graph, 
and a {\it perfect matching\/} is a partial matching 
with the property that every vertex of the graph 
belongs to exactly one edge in the matching.
Here is a perfect matching of
an Aztec diamond graph of order 3:
$$
\pspicture(-1.8,-1.8)(1.8,1.8)
\psset{xunit=.5cm,yunit=.5cm}
\psline(-3,2)(-2,3)
\psline(-3,-2)(-2,-3)
\psline(3,2)(2,3)
\psline(3,-2)(2,-3)
\psline(-3,0)(-2,1)
\psline(0,-3)(1,-2)
\psline(-1,-2)(-2,-1)
\psline(0,-1)(-1,0)
\psline(-1,2)(0,3)
\psline(0,1)(1,2)
\psline(2,-1)(1,0)
\psline(3,0)(2,1)
\endpspicture
$$
(Hereafter, the term ``matching'', used without a modifier,
will denote a perfect matching.)

Matchings of Aztec diamond graphs
were studied in \cite{EKLP},
in the dual guise of domino-tilings
(see below for the domino-tiling
dual to the previously depicted matching).
$$
\pspicture(-2.2,-2.2)(2.2,2.2)
\psset{xunit=.5cm,yunit=.5cm}
\psline(-4,2)(-2,4)
\psline(-4,0)(0,4)
\psline(-4,-2)(2,4)
\psline(0,0)(2,2)
\psline(-2,-4)(1,-1)
\psline(3,1)(4,2)
\psline(0,-4)(4,0)
\psline(2,-4)(4,-2)
\psline(-4,-2)(-2,-4)
\psline(-4,0)(0,-4)
\psline(-4,2)(-3,1)
\psline(-2,0)(0,-2)
\psline(1,-3)(2,-4)
\psline(-2,2)(2,-2)
\psline(-2,4)(-1,3)
\psline(1,1)(4,-2)
\psline(0,4)(4,0)
\psline(2,4)(4,2)
\endpspicture
$$
Note that this is rotated by 45 degrees from
the original picture presented in \cite{EKLP}.
That article showed, by means of four different proofs,
that the number of 
matchings of an Aztec diamond graph of order $n$ 
is $2^{n(n+1)/2}$
(as had been conjectured a decade earlier
in the physics literature \cite{GCZ}).
One of the proofs of this result 
(developed by authors Kuperberg and Propp)
used a geometric-combinatorial procedure dubbed ``shuffling" 
(short for ``domino-shuffling"). 
Shuffling was originally introduced purely 
for purposes of counting
the domino-tilings of the Aztec diamond of order $n$,
but it later turned out to be useful as the basis 
for an algorithm for sampling from 
the uniform distribution on this set of tilings. 
An undergraduate, Sameera Iyengar, 
was the first person to implement
domino-shuffling on a computer,
and the output of her program suggested that 
domino-tilings of Aztec diamonds
exhibit a spatially-expressed phase transition,
where the boundary 
between ``frozen'' and ``non-frozen'' regions 
in a random domino-tiling of a large Aztec diamond
is roughly circular in shape.
A rigorous analysis of the behavior of the shuffling algorithm
led to the first proof (by Jockusch, Propp, and Shor \cite{JPS})
of the asymptotic circularity
of the boundary of the frozen region
(the ``arctic circle theorem'').

Meanwhile, Alexandru Ionescu,
also an undergraduate at the time,
discovered a recurrence relation related to shuffling
that permits one to calculate fairly efficiently, 
for any edge $e$ in an (unweighted) Aztec diamond graph,
the probability that a random matching of the graph
contains the edge $e$.
This theorem allowed Gessel, Ionescu, and Propp
(in unpublished work)
to prove a conjecture of Jockusch
concerning the value of this probability
when $e$ is one of the four central edges in the graph.
More significantly, the theorem allowed 
Cohn, Elkies, and Propp \cite{CEP}
to give a detailed analysis of 
the asymptotic behavior of this probability
as a function of the position of the edge $e$
as the size of the Aztec diamond graph goes to infinity.

During this same period, 
it had become clear that (weighted) enumeration
of matchings of weighted Aztec diamond graphs
had relevance to tiling-models other than domino-tiling.
Notably, Bo-Yin Yang had shown \cite{Y}
that the number of ``diabolo-tilings'' of a ``fortress of order $n$''
was given by a certain formula conjectured by Propp
(for more on this, see section 8).
Yang's proof made use of Kuperberg's re-formulation of the problem
as a question concerning weighted enumeration of matchings
of a weighted version of the Aztec diamond graph
in which some edges had weight 1
while other edges had weight $\frac12$
(see the end of section 5 of this article
and the longer discussion in section 8).
It was clearly desirable to try to extend
the shuffling algorithm to the context of weighted matchings,
or as physicists would call it, the dimer model
in the presence of non-uniform bond-weights,
but at the time it was unclear how to devise the right extension.
Ciucu's work \cite{Ci2}
went part of the way towards this,
and in particular it gave a very clean combinatorial proof
of Yang's result.
However, it was not clear whether one could
efficiently generate diabolo-tilings
using some extension of domino-shuffling.

The purpose of the present article is to give the details
of just such a generalization of domino-shuffling 
into the general context
of Aztec diamond graphs with weighted edges.
At the same time, my goal is to remove some of the mystery
that hung about the shuffling algorithm in the original exposition
in \cite{EKLP},
and to give a exposition of the work of Ionescu,
which up until now has not been described in the literature
at any useful level of detail.
It should be stressed 
\begin{itemize}
\item[$\cdot$]
that the counting algorithm
(presented in section 2 and analyzed in section 5)
is essentially a restatement of Mihai Ciucu's 
cellular graph complementation algorithm \cite{Ci2},
\item[$\cdot$]
that the probability computation algorithm
(presented in section 3 and analyzed in section 6)
is a reasonably straightforward generalization
of Ionescu's unpublished algorithm for
the unweighted case,
\item[$\cdot$]
and that the random generation algorithm
(presented in section 4 and analyzed in section 7),
although algebraically more complicated,
represents no conceptual advance over
the unweighted version of shuffling presented
by Kuperberg and Propp.
\end{itemize}
I think that the major contribution offered by this paper
is not any one of these algorithms in isolation
but the way they fit together in a uniform framework.

The essential tool that makes generalized shuffling work
is a principle first communicated to me by Kuperberg,
which is very similar to tricks common in
the literature on exactly solvable lattice models.
This principle is a lemma that asserts a relationship between
the dimer model on one graph
and the dimer model on another,
where the two graphs differ only in that
a small patch of one graph (sometimes called a ``city'')
has been slightly modified in a particular way
(a process sometimes called ``urban renewal'').
Urban renewal is a powerful trick,
especially when used in combination with an even more trivial trick
called vertex splitting/merging.
(Such graph-rewriting rules are not new to the subject
of enumeration of matchings;
see for instance the article of 
Colbourn, Provan, and Vertigan \cite{CPV}
for a very general approach that uses the ``wye-delta'' lemma
in place of the ``urban renewal'' lemma.)
Urban renewal, applied in $n^2$ locations
in an Aztec diamond of order $n$,
essentially converts it into an Aztec diamond of order $n-1$.
In this way, a relation is established
between matchings of $A_n$ and matchings of $A_{n-1}$.
In particular, one can reduce the weighted enumeration 
of matchings of a weighted Aztec diamond graph of order $n$ 
to a (differently) weighted enumeration of matchings
of the Aztec diamond graph of order $n-1$.

The results described in this article
were first communicated by means of
several long email messages I sent out 
to a number of colleagues in 1996.
I am greatly indebted to Greg Kuperberg,
who in the Fall of 2000 turned these messages
into the first draft of this article,
and who in particular created most of the illustrations
(thereby introducing me to the joys of {\tt{pstricks}}).
Also, Harald Helfgott's program {\tt ren}
(currently available from {\tt http://www.math.wisc.edu/$\sim$propp/};
see the files {\tt ren.c}, {\tt ren.h}, and {\tt ren.html})
was the first implementation of the general algorithm,
and it facilitated several useful discoveries,
including phenomena associated with fortresses \cite{Ke2}.
Helfgott and Ionescu and Iyengar were three 
among many undergraduate research assistants
whose participation has contributed to my work on tilings;
the other members of the Tilings Research Group 
(most of whom were undergraduates at the time
the research was done)
were Pramod Achar, Karen Acquista, Josie Ammer, Federico Ardila,
Rob Blau, Matt Blum, Ruth Britto-Pacumio, Constantin Chiscanu,
Henry Cohn, Chris Douglas, Edward Early, Marisa Gioioso,
David Gupta, Sharon Hollander, 
Julia Khodor, Neelakantan Krishnaswami,
Eric Kuo, Ching Law, Andrew Menard, Alyce Moy, Ben Raphael,
Vis Taraz, Jordan Weitz, Ben Wieland, Lauren Williams,
David Wilson, Jessica Wong, Jason Woolever, and Laurence Yogman. 
I thank Chen Zeng,
who has made use of the generalized shuffling algorithm \cite{ZLH}
in his own work on
and thereby encouraged me to write up this algorithm
for publication.
Lastly, I thank the anonymous referee,
whose careful reading and helpful comments
have made this a better article in every way.

\subsection{Two examples}

Here are two illustrations
of the flexibility of the procedures
described in this article.
The figure below shows a subgraph of
the Aztec diamond graph $A_5$.
This can be thought of as an edge-weighting
of the Aztec diamond graph of order 5
in which the marked edges are assigned weight 1
and the remaining (absent) edges are assigned weight 0.
$$
\pspicture(-1.8,-1.8)(1.8,1.8)
\psset{xunit=.3cm,yunit=.3cm}
\psline(-5,0)(0,5)
\psline(-4,-1)(1,4)
\psline(-3,-2)(2,3)
\psline(-2,-3)(3,2)
\psline(-1,-4)(4,1)
\psline(0,-5)(5,0)
\psline(-5,0)(0,-5)
\psline(-4,1)(1,-4)
\psline(-3,2)(2,-3)
\psline(-2,3)(3,-2)
\psline(-1,4)(4,-1)
\psline(0,5)(5,0)
\psline(-5,-4)(-4,-5)
\psline(-5,-2)(-4,-3)
\psline(-3,-4)(-2,-5)
\psline(5,-4)(4,-5)
\psline(5,-2)(4,-3)
\psline(3,-4)(2,-5)
\psline(-5,4)(-4,5)
\psline(-5,2)(-4,3)
\psline(-3,4)(-2,5)
\psline(5,4)(4,5)
\psline(5,2)(4,3)
\psline(3,4)(2,5)
\endpspicture
$$
It is easy to see that 
the twelve isolated edges all are ``forced'',
in the sense that every matching of the graph
must contain them.
When these forced edges are pruned from the graph,
what remains is the 6-by-6 square grid.
More generally, for any $n$,
one can assign a $0,1$-weighting to the edges
of the Aztec diamond graph of order $2n-1$
so that the matchings of positive weight (i.e., weight 1)
correspond to the matchings of the $2n$-by-$2n$ square grid.
One can use this correspondence in order 
to count matchings of the square grid,
determine inclusion-probabilities for individual edges,
and generate random matchings
(for a technical caveat, see section 9, item 1).

The second example concerns matchings of graphs called
``hexagon honeycomb graphs''.
These graphs have been studied by chemists
as examples of generalized benzene-like hydrocarbons \cite{CG},
and they also
have connections to the study of plane partitions \cite{CLP}.
To construct such a graph, start with a hexagon
with all internal angles measuring 120 degrees
and with sides of respective lengths $a,b,c,a,b,c$,
where $a$, $b$, and $c$ are non-negative integers,
and divide it (in the unique way) into
unit equilateral triangles,
as shown below with $a=b=c=2$.
$$
\pspicture(-1.8,-1.8)(1.8,1.8)
\psset{xunit=.3cm,yunit=.3cm}
\psline(-\fourx,0)(-\twox,6)
\psline(-\threex,-3)(0,6)
\psline(-\twox,-6)(\twox,6)
\psline(0,-6)(\threex,3)
\psline(\twox,-6)(\fourx,0)
\psline(-\fourx,0)(-\twox,-6)
\psline(-\threex,3)(0,-6)
\psline(-\twox,6)(\twox,-6)
\psline(0,6)(\threex,-3)
\psline(\twox,6)(\fourx,0)
\psline(-\twox,6)(\twox,6)
\psline(-\threex,3)(\threex,3)
\psline(-\fourx,0)(\fourx,0)
\psline(-\threex,-3)(\threex,-3)
\psline(-\twox,-6)(\twox,-6)
\endpspicture
$$
Let $V$ be the set of triangles,
and let $E$ be the set of pairs of triangles
sharing an edge.
(Pictorially, we may represent the elements of $V$
by dots centered in the respective triangles,
and the elements of $E$ by segments joining
elements of $V$ at unit distance.)
The graph $(V,E)$ is the hexagon honeycomb graph
$H_{a,b,c}$.
Here, for instance, is the
hexagon honeycomb graph $H_{2,2,2}$:
$$
\pspicture(-1.8,-1.8)(1.8,1.8)
\psset{xunit=.3cm,yunit=.3cm}
\psline(-\threex,-1)(-\threex,1)
\psline(-\twox,-4)(-\twox,-2)
\psline(-\twox,2)(-\twox,4)
\psline(-\onex,-1)(-\onex,1)
\psline(0,-4)(0,-2)
\psline(0,2)(0,4)
\psline(\onex,-1)(\onex,1)
\psline(\twox,-4)(\twox,-2)
\psline(\twox,2)(\twox,4)
\psline(\threex,-1)(\threex,1)
\psline(-\twox,4)(-\onex,5)
\psline(-\threex,1)(-\twox,2)
\psline(0,4)(\onex,5)
\psline(-\onex,1)(0,2)
\psline(-\twox,-2)(-\onex,-1)
\psline(\onex,1)(\twox,2)
\psline(0,-2)(\onex,-1)
\psline(-\onex,-5)(0,-4)
\psline(\twox,-2)(\threex,-1)
\psline(\onex,-5)(\twox,-4)
\psline(\twox,4)(\onex,5)
\psline(\threex,1)(\twox,2)
\psline(0,4)(-\onex,5)
\psline(\onex,1)(0,2)
\psline(\twox,-2)(\onex,-1)
\psline(-\onex,1)(-\twox,2)
\psline(0,-2)(-\onex,-1)
\psline(\onex,-5)(0,-4)
\psline(-\twox,-2)(-\threex,-1)
\psline(-\onex,-5)(-\twox,-4)
\endpspicture
$$
Perfect matchings of such a graph
correspond to tilings of the $a,b,c,a,b,c$ hexagon
by unit rhombuses with angles of 60 and 120 degrees,
each composed of two of the unit equilateral triangles. 

The following figure shows, illustratively,
how $H_{2,2,2}$
can be embedded in the Aztec diamond graph $A_5$:
$$
\pspicture(-1.8,-1.8)(1.8,1.8)
\psset{xunit=.3cm,yunit=.3cm}
\psline(-5,4)(-4,5)
\psline(-3,4)(-2,5)
\psline(-1,4)(0,5)
\psline(0,3)(1,4)
\psline(2,5)(3,4)
\psline(4,5)(5,4)
\psline(-5,2)(-4,3)
\psline(-3,2)(-2,3)
\psline(2,3)(3,2)
\psline(4,3)(5,2)
\psline(-5,0)(-4,1)
\psline(3,0)(4,1)
\psline(4,-1)(5,0)
\psline(2,-3)(3,-2)
\psline(4,-3)(5,-2)
\psline(0,-5)(1,-4)
\psline(2,-5)(3,-4)
\psline(4,-5)(5,-4)
\psline(-5,-2)(-4,-3)
\psline(-3,0)(-2,-1)
\psline(-1,2)(0,1)
\psline(-5,-4)(-4,-5)
\psline(-3,-2)(-2,-3)
\psline(-1,0)(0,-1)
\psline(1,2)(2,1)
\psline(-3,-4)(-2,-5)
\psline(-1,-2)(0,-3)
\psline(1,0)(2,-1)
\psline(-5,-2)(-1,2)
\psline(-5,-4)(1,2)
\psline(-4,-5)(2,1)
\psline(-2,-5)(2,-1)
\endpspicture
$$
More generally, for any $a,b,c$,
one can assign a $0,1$-weighting to the edges
of a suitably large Aztec diamond graph
so that the matchings of positive weight (i.e., weight 1)
correspond to the matchings of $H_{a,b,c}$.
One can use this correspondence in order 
to count matchings of the honeycomb graph,
determine inclusion-probabilities for individual edges,
and generate random matchings
(once again subject to the technical caveat
discussed at the beginning of section 9).

\subsection{Evaluation and comparison of the algorithms}

I have not done a rigorous study of the running-times
of these algorithms under realistic assumptions,
but if one makes the (unrealistic) assumption
that arithmetic operations take constant time
(regardless of the complexity of the numbers involved),
then it is easy to see that the most straightforward
implementation of generalized shuffling,
when applied in an Aztec diamond of order $n$,
takes time roughly $n^3$
(ignoring factors of $\log n$).
We have found that in practice, shuffling is extremely efficient.
The physicist Zeng
has applied the algorithm to the study of the dimer model
in the presence of random edge-weights \cite{ZLH}. 

Alternative methods of enumerating matchings of graphs
are legion;
the simplest method that applies to all planar graphs
is the Pfaffian method of Kasteleyn \cite{Ka},
and there is a variant due to Percus \cite{Pe}
(the permanent-determinant method)
that applies when the graph is bipartite.
For the special case of 
enumerating matchings of Aztec diamond graphs,
there are close to a dozen different
(or at least different-looking) proofs
in the literature.
One that is particularly worthy of mention
is the recent condensation algorithm of Eric Kuo \cite{Kuo}.
Kuo's algorithm appears to be closely related
to the algorithm presented here,
although I have not worked out the details
of the correspondence.
It is also worth remarking upon
the resemblance between Kuo's basic bilinear relation
and the bilinear relation that appears
in section 6
of this article;
both relations are proved in similar ways. 
    
Once one knows how to count matchings of graphs
(or, in the weighted case, to sum the weights
of all the matchings),
there is a trivial way to use this to compute
edge-inclusion probabilities:
apply the counting-algorithm to both the original graph
and the graph from which a selected edge
has been deleted (along with its two endpoints).
The ratio of these (weighted) counts
gives the desired probability.
If one wants to compute only a single edge-inclusion probability,
this method is quite efficient,
and if one uses Kasteleyn's approach,
one need only calculate two Pfaffians;
however, if one wants to compute all edge-inclusion probabilities,
it is wasteful to compute all the Pfaffians
independently of one another,
since the different subgraphs are all very similar.
Wilson \cite{W}
shows how one can organize the calculation more efficiently,
exploiting the similarities between the subgraphs.
I believe that Wilson's algorithm and the approach given here
are likely to be roughly equivalent
in terms of computational difficulty
and numerical stability.
The new approach has the virtue
of being much simpler to code.

Once one knows how to compute edge-inclusion probabilities,
generating a random matching is not hard:
one can cycle through the set of edges,
making decisions about whether to include or not include an edge
in accordance with the associated inclusion-probability
(whose value in general depends on the outcome
of earlier decisions).
The na{}ive way of implementing this,
like the na{}ive way of computing all the inclusion probabilities
in parallel,
requires many determinant-evaluations
or Pfaffian-evaluations
of highly similar matrices;
as in the case of the previous paragraph,
Wilson \cite{W} has shown how to use this
to increase efficiency, resulting in
an $O(N^{1.5})$ algorithm,
where $N$ is the number of vertices.
Note that in the case of Aztec diamonds,
$N=O(n^2)$ and $N^{1.5}=O(n^3)$,
so this is the same approximate running-time
as the algorithm given in this article.
No algorithm faster than Wilson's is known.

\subsection{Overview of article}

The layout of the rest of the article is as follows.
Sections 2 through 7 correspond to the original six installments
of the 1996 e-mail version of the article.
Sections 2 through 4 give algorithms for
weighted enumeration, computation of probabilities,
and random generation;
sections 5 through 7 give the proofs that
these algorithms are valid.
The expository strategy is to use examples wherever possible,
rather than resort to descriptions of the general situation,
especially where this might entail cumbersome notation.
Section 8 applies the methods of the article
to the study of diabolo-tilings of fortresses.
Section 9 concludes the article with a summary
and some open problems.

\section{Computing weight-sums: the algorithm}

Here is the rule (to be immediately followed by an example)
for finding the sum of the weights of all the matchings
of an Aztec diamond graph: 
Given an Aztec diamond graph of order $n$ 
whose edges are marked with weights, 
first decompose the graph into $n^2$ 4-cycles 
(to be called ``cells'' following Mihai Ciucu's coinage \cite{Ci1}), 
and then replace each marked cell
$$
\pspicture(-1,-1)(1,1)
\pspolygon(-1,0)(0,1)(1,0)(0,-1)
\rput[br](-.55, .55){$w$} \rput[bl](.55, .55){$x$}
\rput[tr](-.55,-.55){$y$} \rput[tl](.55,-.55){$z$}
\endpspicture
$$
by the marked cell
$$
\pspicture(-1,-1)(1,1)
\pspolygon(-1,0)(0,1)(1,0)(0,-1)
\rput[br](-.55, .55){$z/(wz+xy)$} \rput[bl](.55, .55){$y/(wz+xy)$}   
\rput[tr](-.55,-.55){$x/(wz+xy)$} \rput[tl](.55,-.55){$w/(wz+xy)$}
\endpspicture
$$
and strip away the outer flank of edges, 
so that an edge-marked Aztec diamond graph of order $n-1$ remains.  
Then (as will be shown in section 5)
the sum of the weights of the matchings of the graph of order $n$ 
equals the sum of the weights of the matchings 
of the graph of order $n-1$ 
times the product of the $n^2$
(different) factors of the form $wz+xy$
(``cell-factors'').  
So, if you perform the reduction process $n$ times, 
obtaining in the end the Aztec diamond of order $0$ 
(for which the sum of the weights of the matchings is $1$), 
you will find that the weighted sum of the matchings of 
the original graph of order $n$ 
is just the product of
$$n^2 + (n-1)^2 + \ldots + 2^2 + 1^2$$
terms, read off from the cells.

Let us try this with the weighted Aztec diamond graph
discussed earlier, in connection with the $4 \times 4$ grid:
$$
\pspicture(-3,-3)(3,3)
\pspolygon(-3,-2)(2,3)(3,2)(-2,-3)
\pspolygon(-3,0)(0,3)(3,0)(0,-3)
\pspolygon(-3,2)(-2,3)(3,-2)(2,-3)
\rput[tl](-2.45,2.45){1} \rput[tr](-1.55,2.45){0}
\rput[tl](-0.45,2.45){1} \rput[tr]( 0.45,2.45){1}
\rput[tl]( 1.55,2.45){0} \rput[tr]( 2.45,2.45){1}
\rput[bl](-2.45,1.55){0} \rput[br](-1.55,1.55){1}
\rput[bl](-0.45,1.55){1} \rput[br]( 0.45,1.55){1}
\rput[bl]( 1.55,1.55){1} \rput[br]( 2.45,1.55){0}
\rput[tl](-2.45,0.45){1} \rput[tr](-1.55,0.45){1}
\rput[tl](-0.45,0.45){1} \rput[tr]( 0.45,0.45){1}
\rput[tl]( 1.55,0.45){1} \rput[tr]( 2.45,0.45){1}
\rput[bl](-2.45,-0.45){1} \rput[br](-1.55,-0.45){1}
\rput[bl](-0.45,-0.45){1} \rput[br]( 0.45,-0.45){1}
\rput[bl]( 1.55,-0.45){1} \rput[br]( 2.45,-0.45){1}
\rput[tl](-2.45,-1.55){0} \rput[tr](-1.55,-1.55){1}
\rput[tl](-0.45,-1.55){1} \rput[tr]( 0.45,-1.55){1}
\rput[tl]( 1.55,-1.55){1} \rput[tr]( 2.45,-1.55){0}
\rput[bl](-2.45,-2.45){1} \rput[br](-1.55,-2.45){0}
\rput[bl](-0.45,-2.45){1} \rput[br]( 0.45,-2.45){1}
\rput[bl]( 1.55,-2.45){0} \rput[br]( 2.45,-2.45){1}
\endpspicture
$$
First we apply the replacement rule described above:
$$
\pspicture(-3,-3)(3,3)
\pspolygon(-3,-2)(2,3)(3,2)(-2,-3)
\pspolygon(-3,0)(0,3)(3,0)(0,-3)
\pspolygon(-3,2)(-2,3)(3,-2)(2,-3)
\rput[tl](-2.45,2.45){1} \rput[tr](-1.55,2.45){0}
\rput[tl](-0.45,2.45){\half} \rput[tr]( 0.45,2.45){\half}
\rput[tl]( 1.55,2.45){0} \rput[tr]( 2.45,2.45){1}
\rput[bl](-2.45,1.55){0} \rput[br](-1.55,1.55){1}
\rput[bl](-0.45,1.55){\half} \rput[br]( 0.45,1.55){\half}
\rput[bl]( 1.55,1.55){1} \rput[br]( 2.45,1.55){0}
\rput[tl](-2.45,0.45){\half} \rput[tr](-1.55,0.45){\half}
\rput[tl](-0.45,0.45){\half} \rput[tr]( 0.45,0.45){\half}
\rput[tl]( 1.55,0.45){\half} \rput[tr]( 2.45,0.45){\half}
\rput[bl](-2.45,-0.45){\half} \rput[br](-1.55,-0.45){\half}
\rput[bl](-0.45,-0.45){\half} \rput[br]( 0.45,-0.45){\half}
\rput[bl]( 1.55,-0.45){\half} \rput[br]( 2.45,-0.45){\half}
\rput[tl](-2.45,-1.55){0} \rput[tr](-1.55,-1.55){1}
\rput[tl](-0.45,-1.55){\half} \rput[tr]( 0.45,-1.55){\half}
\rput[tl]( 1.55,-1.55){1} \rput[tr]( 2.45,-1.55){0}
\rput[bl](-2.45,-2.45){1} \rput[br](-1.55,-2.45){0}
\rput[bl](-0.45,-2.45){\half} \rput[br]( 0.45,-2.45){\half}
\rput[bl]( 1.55,-2.45){0} \rput[br]( 2.45,-2.45){1}
\endpspicture
$$
Next we strip away the outer layer (note that some of the values
computed at the previous step are simply thrown away, so that
it is not in fact strictly necessary to have computed them):
$$
\pspicture(-2,-2)(2,2)
\pspolygon(-2,-1)(1,2)(2,1)(-1,-2)
\pspolygon(-2,1)(-1,2)(2,-1)(1,-2)
\rput[tl](-1.45,1.45){1} \rput[tr](-0.55,1.45){\half}
\rput[tl](0.55,1.45){\half} \rput[tr](1.45,1.45){1}
\rput[bl](-1.45,0.55){\half} \rput[br](-0.55,0.55){\half}
\rput[bl](0.55,0.55){\half} \rput[br](1.45,0.55){\half}
\rput[tl](-1.45,-0.55){\half} \rput[tr](-0.55,-0.55){\half}
\rput[tl](0.55,-0.55){\half} \rput[tr](1.45,-0.55){\half}
\rput[bl](-1.45,-1.45){1} \rput[br](-0.55,-1.45){\half}
\rput[bl](0.55,-1.45){\half} \rput[br](1.45,-1.45){1}
\endpspicture
$$
Now we do it again (only this time I will not bother to depict
the weights that end up getting thrown away):
$$
\pspicture(-1,-1)(1,1)
\pspolygon(-1,0)(0,1)(1,0)(0,-1)
\rput[br](-.55, .55){4/3} \rput[bl](.55, .55){4/3}
\rput[tr](-.55,-.55){4/3} \rput[tl](.55,-.55){4/3}
\endpspicture
$$
We do it one more time, and we obtain the empty Aztec diamond.

Now we multiply together all the cell-factors:
$$2^5  \times (3/4)^4 \times (32/9) = 36.$$

A different example arises when 
the cells of an Aztec diamond of order $n$
are alternately colored white and black,
with all the edges in the white cells
being assigned one weight ($a$)
and all the edges in the black cells
being assigned another weight ($b$).
In this case, a single application of the algorithm
gives an Aztec diamond of order $n-1$
in which there are again edges of two different weights,
only now the different-weighted edges
are intermixed in the cells.
However, if one repeats the algorithm again,
one gets an Aztec diamond of order $n-2$
in which the alternation of the two weights
is as in the original graph.
In this fashion, it is possible to express
the sum of the weights of all the matchings of the graph
in terms of cell-factors equal to 
$2a^2$, $2b^2$, and $a^2+b^2$.
Indeed, by specializing these weights with $a=1$ and $b=\frac12$,
we get a very simple proof of the ``powers-of-5'' formula
proved by Yang \cite{Y},
similar to the proof given by Ciucu \cite{Ci2}.

\section{Computing edge-inclusion probabilities: the algorithm}

Let us assume that non-negative real weights have been assigned
to the edges of an Aztec diamond graph,
and let us define the weight of a matching of the graph
as the product of the weights of the constituent edges.
Assume that the weighting of the edges is such that
there exists at least one matching of non-zero weight.
Then the weighting of the matchings
induces a probability distribution 
on the set of matchings of the graph,
in which the probability of a particular matching is proportional
to the weight of that matching.
This is most natural in the case where all edges
have weight 0 or weight 1;
then one is looking at the uniform distribution on
the set of all matchings of the subgraph
consisting of the edges of weight 1.
In any case, the probability distribution is well-defined,
so it makes sense to ask,
what is the probability that a random matching
(chosen relative to the aforementioned probability distribution)
will contain some particular edge $e$?
Here we present a scheme that simultaneously answers this question
for all edges $e$ of the graph.

What we do is thread our way back 
through the reduction process described in the previous section, 
starting from order 0 and working our way back up to order $n$,
computing the edge probability in successively larger and larger 
weighted Aztec diamonds graphs;
along the way we make use of the cell-weights 
that were computed during the reduction process.  
Suppose we have computed edge probabilities for
the weighted graph of order $n-1$.  
To derive the edge probabilities
for the weighted graph of order $n$, 
we embed the smaller graph in
the larger (concentrically), 
divide the larger graph into $n^2$ cells,
and swap the numbers belonging to edges 
on opposite sides of a cell
(the example below will make this clearer).  
When we have done this,
the numbers on the edges are (typically) 
not yet the true probabilities, 
but they can be regarded as approximations to them, 
in the sense that we can add a slightly more complicated
correction term that makes the ``approximate'' formulas exact.  
Let us zoom in on a particular cell to see
how it works: 
Consider the typical cell
$$
\pspicture(-1,-1)(1,1)
\pspolygon(-1,0)(0,1)(1,0)(0,-1)
\rput[br](-.55, .55){$w$} \rput[bl](.55, .55){$x$}
\rput[tr](-.55,-.55){$y$} \rput[tl](.55,-.55){$z$}
\endpspicture
$$
and suppose that the edges with weights $w$, $x$, $y$, and $z$ 
(before the swapping has taken place) 
have been given ``approximate'' probabilities
$p$, $q$, $r$, and $s$, respectively.  
Then the \emph{exact} probabilities for
these respective edges are
\begin{align*}
&s + \frac{(1-p-q-r-s)wz}{wz+xy} && r + \frac{(1-p-q-r-s)xy}{wz+xy} \\
&q + \frac{(1-p-q-r-s)xy}{wz+xy} && p + \frac{(1-p-q-r-s)wz}{wz+xy}.
\end{align*}
The number $1-p-q-r-s$ will be called the \emph{deficit}
associated with that cell, 
and in the context of the four preceding inset expressions,
it is called the \emph{(net) creation rate};
the numbers $wz/(wz+xy)$ and $xy/(wz+xy)$ 
are called \emph{(net) creation biases}.

Applying this to the $4 \times 4$ grid graph, 
we find that for the weighted Aztec graphs 
of orders $1$, $2$, and $3$ that we found (in reverse order)
when we \emph{counted} the matchings, 
the edge-probabilities are as follows:
For the order $1$ graph, we have
$$
\pspicture(-1,-1)(1,1)
\pspolygon(-1,0)(0,1)(1,0)(0,-1)
\rput[br](-.55, .55){\half} \rput[bl](.55, .55){\half}
\rput[tr](-.55,-.55){\half} \rput[tl](.55,-.55){\half}
\endpspicture
$$
We embed this in the Aztec graph of order $2$:
$$
\pspicture(-2,-2)(2,2)
\pspolygon(-2,-1)(1,2)(2,1)(-1,-2)
\pspolygon(-2,1)(-1,2)(2,-1)(1,-2)
\qdisk(-1,-1){.05} \qdisk(-1,1){.05} \qdisk(1,-1){.05} \qdisk(1,1){.05}
\rput[tl](-1.45,1.45){0} \rput[tr](-0.55,1.45){0}
\rput[tl](0.55,1.45){0} \rput[tr](1.45,1.45){0}
\rput[bl](-1.45,0.55){0} \rput[br](-0.55,0.55){\half}
\rput[bl](0.55,0.55){\half} \rput[br](1.45,0.55){0}
\rput[tl](-1.45,-0.55){0} \rput[tr](-0.55,-0.55){\half}
\rput[tl](0.55,-0.55){\half} \rput[tr](1.45,-0.55){0}
\rput[bl](-1.45,-1.45){0} \rput[br](-0.55,-1.45){0}
\rput[bl](0.55,-1.45){0} \rput[br](1.45,-1.45){0}
\endpspicture
$$
Note that I have put a dot in the middle of each of the four cells.
We swap each number with the number opposite it in its cell:
$$
\pspicture(-2,-2)(2,2)
\pspolygon(-2,-1)(1,2)(2,1)(-1,-2)
\pspolygon(-2,1)(-1,2)(2,-1)(1,-2)
\qdisk(-1,-1){.05} \qdisk(-1,1){.05} \qdisk(1,-1){.05} \qdisk(1,1){.05}
\rput[tl](-1.45,1.45){\half} \rput[tr](-0.55,1.45){0}
\rput[tl](0.55,1.45){0} \rput[tr](1.45,1.45){\half}
\rput[bl](-1.45,0.55){0} \rput[br](-0.55,0.55){0}
\rput[bl](0.55,0.55){0} \rput[br](1.45,0.55){0}
\rput[tl](-1.45,-0.55){0} \rput[tr](-0.55,-0.55){0}
\rput[tl](0.55,-0.55){0} \rput[tr](1.45,-0.55){0}
\rput[bl](-1.45,-1.45){\half} \rput[br](-0.55,-1.45){0}
\rput[bl](0.55,-1.45){0} \rput[br](1.45,-1.45){\half}
\endpspicture
$$
These are the approximate, inexact edge-probabilities.  
To find the exact values, we compute that in each of the four cells 
the deficit $1-p-q-r-s$ is $1/2$.  
The weights on the edges 
(all equal to $1$ or $1/2$ --- see second-to-list figure of section 2) 
give us creation biases
$$(1/2)/(3/4) = 2/3$$
and
$$(1/4)/(3/4) = 1/3.$$
So we increment the $1/2$'s (and the $0$'s that are opposite
them) by
$$(1/2)(2/3) = 1/3,$$
and we increment the other $0$'s by
$$(1/2)(1/3) = 1/6,$$
obtaining
$$
\pspicture(-2,-2)(2,2)
\pspolygon(-2,-1)(1,2)(2,1)(-1,-2)
\pspolygon(-2,1)(-1,2)(2,-1)(1,-2)
\rput[tl](-1.45, 1.45){\small 5/6} \rput[tr](-0.55, 1.45){\small 1/6}
\rput[tl]( 0.55, 1.45){\small 1/6} \rput[tr]( 1.45, 1.45){\small 5/6}
\rput[bl](-1.45, 0.55){\small 1/6} \rput[br](-0.55, 0.55){\small 1/3}
\rput[bl]( 0.55, 0.55){\small 1/3} \rput[br]( 1.45, 0.55){\small 1/6}
\rput[tl](-1.45,-0.55){\small 1/6} \rput[tr](-0.55,-0.55){\small 1/3}
\rput[tl]( 0.55,-0.55){\small 1/3} \rput[tr]( 1.45,-0.55){\small 1/6}
\rput[bl](-1.45,-1.45){\small 5/6} \rput[br](-0.55,-1.45){\small 1/6}
\rput[bl]( 0.55,-1.45){\small 1/6} \rput[br]( 1.45,-1.45){\small 5/6}
\endpspicture
$$
Embed this in the graph of order $3$, 
and swap the numbers in each cell:
$$
\pspicture(-3,-3)(3,3)
\pspolygon(-3,-2)(2,3)(3,2)(-2,-3)
\pspolygon(-3,0)(0,3)(3,0)(0,-3)
\pspolygon(-3,2)(-2,3)(3,-2)(2,-3)
\qdisk(-2,-2){.05}\qdisk(-2, 0){.05}\qdisk(-2, 2){.05}
\qdisk( 0,-2){.05}\qdisk( 0, 0){.05}\qdisk( 0, 2){.05}
\qdisk( 2,-2){.05}\qdisk( 2, 0){.05}\qdisk( 2, 2){.05}
\rput[tl](-2.45, 2.45){\small 5/6} \rput[tr](-1.55, 2.45){0}
\rput[tl](-0.45, 2.45){\small 1/6} \rput[tr]( 0.45, 2.45){\small 1/6}
\rput[tl]( 1.55, 2.45){0} \rput[tr]( 2.45, 2.45){\small 5/6}
\rput[bl](-2.45, 1.55){0} \rput[br](-1.55, 1.55){0}
\rput[bl](-0.45, 1.55){0} \rput[br]( 0.45, 1.55){0}
\rput[bl]( 1.55, 1.55){0} \rput[br]( 2.45, 1.55){0}
\rput[tl](-2.45, 0.45){\small 1/6} \rput[tr](-1.55, 0.45){0}
\rput[tl](-0.45, 0.45){\small 1/3} \rput[tr]( 0.45, 0.45){\small 1/3}
\rput[tl]( 1.55, 0.45){0} \rput[tr]( 2.45, 0.45){\small 1/6}
\rput[bl](-2.45,-0.45){\small 1/6} \rput[br](-1.55,-0.45){0}
\rput[bl](-0.45,-0.45){\small 1/3} \rput[br]( 0.45,-0.45){\small 1/3}
\rput[bl]( 1.55,-0.45){0} \rput[br]( 2.45,-0.45){\small 1/6}
\rput[tl](-2.45,-1.55){0} \rput[tr](-1.55,-1.55){0}
\rput[tl](-0.45,-1.55){0} \rput[tr]( 0.45,-1.55){0}
\rput[tl]( 1.55,-1.55){0} \rput[tr]( 2.45,-1.55){0}
\rput[bl](-2.45,-2.45){\small 5/6} \rput[br](-1.55,-2.45){0}
\rput[bl](-0.45,-2.45){\small 1/6} \rput[br]( 0.45,-2.45){\small 1/6}
\rput[bl]( 1.55,-2.45){0} \rput[br]( 2.45,-2.45){\small 5/6}
\endpspicture
$$
The four corner-cells have deficit $1/6$, 
the other four outer cells have deficit $2/3$ 
and the inner cell has deficit $-1/3$ 
(or surplus $+1/3$).  So we do our adjustments, 
this time with all creation biases equal to $1/2$ 
(except in the corners, 
where the creation biases are $0$ and $1$ ---
see third-to-list figure of section 2): 
$$
\pspicture(-3,-3)(3,3)
\pspolygon(-3,-2)(2,3)(3,2)(-2,-3)
\pspolygon(-3,0)(0,3)(3,0)(0,-3)
\pspolygon(-3,2)(-2,3)(3,-2)(2,-3)
\qdisk(-2,-2){.05}\qdisk(-2, 0){.05}\qdisk(-2, 2){.05}
\qdisk( 0,-2){.05}\qdisk( 0, 0){.05}\qdisk( 0, 2){.05}
\qdisk( 2,-2){.05}\qdisk( 2, 0){.05}\qdisk( 2, 2){.05}
\rput[tl](-2.45, 2.45){1} \rput[tr](-1.55, 2.45){0}
\rput[tl](-0.45, 2.45){\small 1/2} \rput[tr]( 0.45, 2.45){\small 1/2}
\rput[tl]( 1.55, 2.45){0} \rput[tr]( 2.45, 2.45){1}
\rput[bl](-2.45, 1.55){0} \rput[br](-1.55, 1.55){\small 1/6}
\rput[bl](-0.45, 1.55){\small 1/3} \rput[br]( 0.45, 1.55){\small 1/3}
\rput[bl]( 1.55, 1.55){\small 1/6} \rput[br]( 2.45, 1.55){0}
\rput[tl](-2.45, 0.45){\small 1/2} \rput[tr](-1.55, 0.45){\small 1/3}
\rput[tl](-0.45, 0.45){\small 1/6} \rput[tr]( 0.45, 0.45){\small 1/6}
\rput[tl]( 1.55, 0.45){\small 1/3} \rput[tr]( 2.45, 0.45){\small 1/2}
\rput[bl](-2.45,-0.45){\small 1/2} \rput[br](-1.55,-0.45){\small 1/3}
\rput[bl](-0.45,-0.45){\small 1/6} \rput[br]( 0.45,-0.45){\small 1/6}
\rput[bl]( 1.55,-0.45){\small 1/3} \rput[br]( 2.45,-0.45){\small 1/2}
\rput[tl](-2.45,-1.55){0} \rput[tr](-1.55,-1.55){\small 1/6}
\rput[tl](-0.45,-1.55){\small 1/3} \rput[tr]( 0.45,-1.55){\small 1/3}
\rput[tl]( 1.55,-1.55){\small 1/6} \rput[tr]( 2.45,-1.55){0}
\rput[bl](-2.45,-2.45){1} \rput[br](-1.55,-2.45){0}
\rput[bl](-0.45,-2.45){\small 1/2} \rput[br]( 0.45,-2.45){\small 1/2}
\rput[bl]( 1.55,-2.45){0} \rput[br]( 2.45,-2.45){1}
\endpspicture
$$
These are the true edge-inclusion probabilities 
associated with the original graph.  
That is, if we take a uniform random matching 
of the $4 \times 4$ grid, 
these are the probabilities with which we will see 
the respective edges occur in the matching.

\section{Generating a random matching: the algorithm}

Before one can begin to generate random matchings of
a weighted Aztec diamond graph of order $n$, 
one must apply the reduction algorithm used in the counting algorithm,
obtaining weighted Aztec diamond graphs of
orders $n-1$, $n-2$, etc.  
But once this has been done, 
generating random matchings of the graph is quite easy.  
The algorithm is an iterative one: 
starting from a Aztec diamond graph of order $0$, 
one successively generates random matchings 
of the graphs of orders $1$, $2$, $3$, etc., 
using the weights that one computed during the reduction-process.  
At each stage, the probability of seeing any particular matching 
can be shown to be proportional to the weight of that matching 
(where the constant of proportionality naturally changes 
as one progresses to larger and larger Aztec diamond graphs).

Here is how one iteration of the procedure works.  
Given a perfect matching of the Aztec diamond graph of order $k-1$, 
embed the matching inside the Aztec diamond graph of order $k$, 
so that you have a partial matching of an Aztec diamond graph order $k$.  
(Note that the $4$-cycles that were cells in the small graph
become the holes between the cells in the large graph).  
In the new enlarged graph, find all (new) cells 
that contain two matched edges and delete both edges. 
(In \cite{EKLP}, this was called ``destruction''.) 
Then replace each edge in the resulting matching 
with the edge opposite it in its cell.  
(In \cite{EKLP}, this was called ``shuffling'',
although in subsequent talks and articles
I prefer to call this step of the algorithm ``sliding''
and to reserve the term ``shuffling'' for
the compound operation consisting
of ``creation'', ``sliding'', and ``destruction''.)
It can be shown that the complement of the resulting matching 
(that is, the graph that remains when the matched vertices
and all their incident edges are removed) 
can be covered by $4$-cycles in a unique way 
(and it is easy to find this covering 
just by making a one-time scan through the graph).  
Each such $4$-cycle has weights attached to each of its $4$ edges, 
so we can choose a random perfect matching of the 4-cycle
using the weights to determine the respective probabilities 
of the two different ways to match the cycle.
(In \cite{EKLP}, this was called ``creation''.)  
Taking the union of these new edges with
the edges that are already in place, 
we get a perfect matching of the Aztec
diamond graph of order $k$.

Let us see how this works with the $4 \times 4$ grid 
(for which we have already worked out the reduction process).  
To generate a random perfect matching, 
we start with the empty matching of the graph of order 0, 
embed it in the graph of order 1, 
and find that the complementary graph of the empty matching 
is a single 4-cycle, which we match in one of two possible ways. 
In this case,
each of the four edges has weight $4/3$, 
so both matchings have equal likelihood; say we choose
$$
\pspicture(-1,-1)(1,1)
\psline(-1,0)(0,1)\psline(0,-1)(1,0)
\qdisk(0,0){.05} \odisk( 1,0){.1} \odisk(-1,0){.1}
\odisk( 0,1){.1} \odisk( 0,-1){.1}
\endpspicture
$$
(Here we have used open circles to denote the vertices of the graph.)
Now we embed this matching in the Aztec diamond graph of order $2$:
$$
\pspicture(-2,-2)(2,2)
\psline(-1,0)(0,1)\psline(0,-1)(1,0)
\qdisk(-1,-1){.05}\qdisk(1,-1){.05}
\qdisk(-1,1){.05}\qdisk(1,1){.05}
\odisk(-1,2){.1} \odisk( 1,2){.1}
\odisk(-2,1){.1} \odisk( 0,1){.1} \odisk( 2,1){.1}
\odisk(-1,0){.1} \odisk( 1,0){.1}
\odisk(-2,-1){.1} \odisk( 0,-1){.1} \odisk( 2,-1){.1}
\odisk(-1,-2){.1} \odisk( 1,-2){.1}
\endpspicture
$$
Neither of the two edges needs to be destroyed, because they belong
to different cells.  Replacing each by the opposite edge in its cell,
we get
$$
\pspicture(-2,-2)(2,2)
\psline(-2,1)(-1,2)\psline(2,-1)(1,-2)
\qdisk(-1,-1){.05}\qdisk(1,-1){.05}
\qdisk(-1,1){.05}\qdisk(1,1){.05}
\odisk(-1,2){.1} \odisk( 1,2){.1}
\odisk(-2,1){.1} \odisk( 0,1){.1} \odisk( 2,1){.1}
\odisk(-1,0){.1} \odisk( 1,0){.1}
\odisk(-2,-1){.1} \odisk( 0,-1){.1} \odisk( 2,-1){.1}
\odisk(-1,-2){.1} \odisk( 1,-2){.1}
\endpspicture
$$
The complementary graph has a unique cover by $4$-cycles:
$$
\pspicture(-2,-2)(2,2)
\pspolygon(-2,-1)(-1,0)(0,-1)(-1,-2)
\pspolygon(2,1)(1,0)(0,1)(1,2)
\qdisk(-1,-1){.05}\qdisk(1,-1){.05}
\qdisk(-1,1){.05}\qdisk(1,1){.05}
\odisk(-1,2){.1} \odisk( 1,2){.1}
\odisk(-2,1){.1} \odisk( 0,1){.1} \odisk( 2,1){.1}
\odisk(-1,0){.1} \odisk( 1,0){.1}
\odisk(-2,-1){.1} \odisk( 0,-1){.1} \odisk( 2,-1){.1}
\odisk(-1,-2){.1} \odisk( 1,-2){.1}
\endpspicture
$$  
In each $4$-cycle, three edges have weight $1/2$ 
and one has weight $1$,
so we are twice as likely to see the one of the matchings as the other.  
Let us suppose that when we choose matchings with suitable bias, 
we get the more likely of the two matchings in both cells 
(as happens $4/9$ of the time).  Thus we have
$$
\pspicture(-2,-2)(2,2)
\qdisk(-1,-1){.05}\qdisk(1,-1){.05}
\qdisk(-1,1){.05}\qdisk(1,1){.05}
\psline(-2,1)(-1,2)\psline(2,1)(1,2)\psline(1,0)(0,1)
\psline(2,-1)(1,-2)\psline(-2,-1)(-1,-2)\psline(-1,0)(0,-1)
\odisk(-1,2){.1} \odisk( 1,2){.1}
\odisk(-2,1){.1} \odisk( 0,1){.1} \odisk( 2,1){.1}
\odisk(-1,0){.1} \odisk( 1,0){.1}
\odisk(-2,-1){.1} \odisk( 0,-1){.1} \odisk( 2,-1){.1}
\odisk(-1,-2){.1} \odisk( 1,-2){.1}
\endpspicture
$$
Now we embed the matching in the graph of order $3$:
$$
\pspicture(-3,-3)(3,3)
\qdisk(-2,2){.05}\qdisk(0,2){.05}\qdisk(2,2){.05}
\qdisk(-2,0){.05}\qdisk(0,0){.05}\qdisk(2,0){.05}
\qdisk(-2,-2){.05}\qdisk(0,-2){.05}\qdisk(2,-2){.05}
\psline(-2,1)(-1,2)\psline(2,1)(1,2)\psline(1,0)(0,1)
\psline(2,-1)(1,-2)\psline(-2,-1)(-1,-2)\psline(-1,0)(0,-1)
\odisk(-2,3){.1} \odisk( 0,3){.1} \odisk( 2,3){.1}
\odisk(-3,2){.1} \odisk(-1,2){.1} \odisk( 1,2){.1} \odisk( 3,2){.1}
\odisk(-2,1){.1} \odisk( 0,1){.1} \odisk( 2,1){.1}
\odisk(-3,0){.1} \odisk(-1,0){.1} \odisk( 1,0){.1} \odisk( 3,0){.1}
\odisk(-2,-1){.1} \odisk( 0,-1){.1} \odisk( 2,-1){.1}
\odisk(-3,-2){.1} \odisk(-1,-2){.1} \odisk( 1,-2){.1} \odisk( 3,-2){.1}
\odisk(-2,-3){.1} \odisk( 0,-3){.1} \odisk( 2,-3){.1}
\endpspicture
$$
This time we must delete the two central edges, 
because they belong to the same cell.  
After we have done this, 
we replace each of the four remaining edges 
with the opposite edge of its $4$-cycle, obtaining
$$
\pspicture(-3,-3)(3,3)
\qdisk(-2,2){.05}\qdisk(0,2){.05}\qdisk(2,2){.05}
\qdisk(-2,0){.05}\qdisk(0,0){.05}\qdisk(2,0){.05}
\qdisk(-2,-2){.05}\qdisk(0,-2){.05}\qdisk(2,-2){.05}
\psline(-3,2)(-2,3)\psline(2,3)(3,2)
\psline(-3,-2)(-2,-3)\psline(2,-3)(3,-2)
\odisk(-2,3){.1} \odisk( 0,3){.1} \odisk( 2,3){.1}
\odisk(-3,2){.1} \odisk(-1,2){.1} \odisk( 1,2){.1} \odisk( 3,2){.1}
\odisk(-2,1){.1} \odisk( 0,1){.1} \odisk( 2,1){.1}
\odisk(-3,0){.1} \odisk(-1,0){.1} \odisk( 1,0){.1} \odisk( 3,0){.1}
\odisk(-2,-1){.1} \odisk( 0,-1){.1} \odisk( 2,-1){.1}
\odisk(-3,-2){.1} \odisk(-1,-2){.1} \odisk( 1,-2){.1} \odisk( 3,-2){.1}
\odisk(-2,-3){.1} \odisk( 0,-3){.1} \odisk( 2,-3){.1}
\endpspicture
$$
There is a unique cover of the complement by $4$-cycles:
$$
\pspicture(-3,-3)(3,3)
\qdisk(-2,2){.05}\qdisk(0,2){.05}\qdisk(2,2){.05}
\qdisk(-2,0){.05}\qdisk(0,0){.05}\qdisk(2,0){.05}
\qdisk(-2,-2){.05}\qdisk(0,-2){.05}\qdisk(2,-2){.05}
\pspolygon(-3,0)(-2,1)(-1,0)(-2,-1)
\pspolygon(3,0)(2,-1)(1,0)(2,1)
\pspolygon(0,-3)(1,-2)(0,-1)(-1,-2)
\pspolygon(0,3)(-1,2)(0,1)(1,2)
\odisk(-2,3){.1} \odisk( 0,3){.1} \odisk( 2,3){.1}
\odisk(-3,2){.1} \odisk(-1,2){.1} \odisk( 1,2){.1} \odisk( 3,2){.1}
\odisk(-2,1){.1} \odisk( 0,1){.1} \odisk( 2,1){.1}
\odisk(-3,0){.1} \odisk(-1,0){.1} \odisk( 1,0){.1} \odisk( 3,0){.1}
\odisk(-2,-1){.1} \odisk( 0,-1){.1} \odisk( 2,-1){.1}
\odisk(-3,-2){.1} \odisk(-1,-2){.1} \odisk( 1,-2){.1} \odisk( 3,-2){.1}
\odisk(-2,-3){.1} \odisk( 0,-3){.1} \odisk( 2,-3){.1}
\endpspicture
$$
As it happens, each edge in the complement has weight $1$ 
in the (original) weighted Aztec diamond of order $3$, 
so we can use a fair coin to decide 
how to match each of the $4$-cycles.

\section{Computing weight-sums: the proof}

To prove the claim, 
we begin with a lemma 
(the ``urban renewal lemma"):
If you have a weighted graph $G$
that contains the local pattern
$$
\pspicture(-2.5,-2.5)(2.5,2.5)
\pspolygon(-1,0)(0,1)(1,0)(0,-1)
\psline(1,0)(2.4,0)\psline(-1,0)(-2.4,0)
\psline(0,1)(0,2.4)\psline(0,-1)(0,-2.4)
\psline(.2,2.4)(0,2)(-.2,2.4)
\psline(.2,-2.4)(0,-2)(-.2,-2.4)
\psline(2.4,.2)(2,0)(2.4,-.2)
\psline(-2.4,.2)(-2,0)(-2.4,-.2)
\rput[br](-.55, .55){$w$} \rput[bl](.55, .55){$x$}
\rput[tr](-.55,-.55){$y$} \rput[tl](.55,-.55){$z$}
\odisk(0,-2){.1} \odisk(0,-1){.1} \odisk(0, 1){.1} \odisk(0, 2){.1}
\odisk(-2,0){.1} \odisk(-1,0){.1} \odisk( 1,0){.1} \odisk( 2,0){.1}
\rput[r](-.2,2){$A$}
\rput[b](-2,.2){$B$}
\rput[b](2,.2){$C$}
\rput[r](-.2,-2){$D$}
\endpspicture
$$
(here $A$,$B$,$C$,$D$ are vertex labels 
and $w$,$x$,$y$,$z$ are edge weights,
with unmarked edges having weight $1$) and you define $G'$ as the
graph you get when you replace this local pattern by
$$
\pspicture(-2.5,-2.5)(2.5,2.5)
\pspolygon(-2,0)(0,2)(2,0)(0,-2)
\psline(2,0)(2.4,0)\psline(-2,0)(-2.4,0)
\psline(0,2)(0,2.4)\psline(0,-2)(0,-2.4)
\psline(.2,2.4)(0,2)(-.2,2.4)
\psline(.2,-2.4)(0,-2)(-.2,-2.4)
\psline(2.4,.2)(2,0)(2.4,-.2)
\psline(-2.4,.2)(-2,0)(-2.4,-.2)
\rput[br](-1.05, 1.05){$w'$} \rput[bl](1.05, 1.05){$x'$}
\rput[tr](-1.05,-1.05){$y'$} \rput[tl](1.05,-1.05){$z'$}
\odisk(0,-2){.1} \odisk(0,2){.1}
\odisk(-2,0){.1} \odisk(2,0){.1}
\rput[r](-.2,2){$A$}
\rput[b](-2,.2){$B$}
\rput[b](2,.2){$C$}
\rput[r](-.2,-2){$D$}
\endpspicture
$$
with
\begin{align*}
w' &= z/(wz+xy) & x' &= y/(wz+xy) \\
y' &= x/(wz+xy) & z' &= w/(wz+xy),
\end{align*}
then the sum of the weights of the matchings of $G$ equals $wz+xy$ times
the sum of the weights of the matchings of $G'$.  The proof will be
by verification:
For each possible subset $S$ of $\{A,B,C,D\}$,
    we will check that
	the total weight of the matchings of $G$ in which
	    the vertices in $S$ are matched with vertices 
    inside the patch
	    and the vertices in $\{A,B,C,D\}\setminus S$ are matched 
    outside the patch
	equals $wz+xy$ times the total weight of the matchings of $G'$ 
    in which the vertices in $S$ are matched with vertices 
    inside the patch
	    and the vertices in $\{A,B,C,D\}\setminus S$ are matched 
    outside the patch.
As it turns out, $10$ of the $2^4=16$ cases
are trivial, and of the remaining $6$, 
four are related by symmetry, 
so it comes down to three cases.

\begin{description}
\item[$\bullet$] $S = \{A,B\}$: 
In this case our matching of $G$ 
matches $A$,$B$ inward and $C$,$D$ outward, 
so it must contain the edge of weight $z$.
This matching of $G$ is associated with a matching of $G'$ 
in which $A$ and $B$ are matched (inward) to each other
via an edge of weight $w'$
and $C$ and $D$ are matched outward just as before.
The ratio of the weight of the (given) matching of $G$
to the (associated) matching of $G'$ is $z/w' = wz+xy$.

\item[$\bullet$] $S = \{A,B,C,D\}$: 
In this case our matching of $G$
matches all four vertices $A,B,C,D$ inward.
This matching is associated with two matchings of $G'$,
one of which matches $A$ with $B$ and $C$ with $D$
and the other of which matches $A$ with $C$ and $B$ with $D$
(all edges outside of the patch remain as before).
If we let $Q$ denote the weight 
of the given matching of $G$,
one of the two derived matchings of $G'$ has weight $w'z'Q$
and the other has weight $x'y'Q$.  
The combined weight of the two matchings is
$$(w'z'+x'y')Q = Q/(wz+xy),$$
so the ratio of the weight of the single (given) matching of $G$
to the two (associated) matchings of $G'$ is $wz+xy$.

\item[$\bullet$] $S = \emptyset$: 
In this case our matching of $G$
matches all four vertices $A,B,C,D$ outward,
and either contains edges $w$ and $z$
or contains edges $x$ and $y$.
These two cases must be lumped together.
Let $Q$ denote the product of the edges
of all the weights of the other edges of the matching
(not counting $w,x,y,z$).
We find that two matchings of $G$, of weight $wzQ$ and $xyQ$, 
are associated with a single matching of $G'$ of weight $Q$.  
The ratio of the weights of the two matchings of $G$
to the single matching of $G'$ is $wz+xy$.
\end{description}

This completes the proof of the lemma.

Now, to prove the Aztec reduction theorem, suppose we have a weighted
Aztec graph of order $n$.
$$
\pspicture(-1.5,-1.5)(1.5,1.5)
\pspolygon(-1.5,-1)(-1,-1.5)(1.5,1)(1,1.5)
\pspolygon(-1.5,0)(0,-1.5)(1.5,0)(0,1.5)
\pspolygon(-1.5,1)(1,-1.5)(1.5,-1)(-1,1.5)
\odisk(-1,1.5){.1}\odisk(0,1.5){.1}\odisk(1,1.5){.1}
\odisk(-1.5,1){.1}\odisk(-.5,1){.1}\odisk(.5,1){.1}\odisk(1.5,1){.1}
\odisk(-1,.5){.1}\odisk(0,.5){.1}\odisk(1,.5){.1}
\odisk(-1.5,0){.1}\odisk(-.5,0){.1}\odisk(.5,0){.1}\odisk(1.5,0){.1}
\odisk(-1,-.5){.1}\odisk(0,-.5){.1}\odisk(1,-.5){.1}
\odisk(-1.5,-1){.1}\odisk(-.5,-1){.1}\odisk(.5,-1){.1}\odisk(1.5,-1){.1}
\odisk(-1,-1.5){.1}\odisk(0,-1.5){.1}\odisk(1,-1.5){.1}
\endpspicture
$$
(Here $n=3$ and the weights are not shown.)
We begin our reduction by splitting each vertex into three
vertices:
$$
\pspicture(-3.5,-3.5)(3.5,3.5)
\psline(-2,2.5)(-2,3.5)\psline(0,2.5)(0,3.5)\psline(2,2.5)(2,3.5)
\psline(2.5,-2)(3.5,-2)\psline(2.5,0)(3.5,0)\psline(2.5,2)(3.5,2)
\multips(-2,-2)(2,0){3}{\multips(0,0)(0,2){3}{
    \psline(-1.5,0)(-.5,0)\psline(0,-1.5)(0,-.5)
    \pspolygon(0,.5)(.5,0)(0,-.5)(-.5,0)
    \odisk(0,.5){.1}\odisk(.5,0){.1}\odisk(0,-.5){.1}\odisk(-.5,0){.1}
    \odisk(0,-1){.1}\odisk(0,-1.5){.1}\odisk(-1,0){.1}\odisk(-1.5,0){.1}}}
\odisk(-2,3){.1}\odisk(0,3){.1}\odisk(2,3){.1}
\odisk(-2,3.5){.1}\odisk(0,3.5){.1}\odisk(2,3.5){.1}
\odisk(3,-2){.1}\odisk(3,0){.1}\odisk(3,2){.1}
\odisk(3.5,-2){.1}\odisk(3.5,0){.1}\odisk(3.5,2){.1}
\endpspicture
$$
It is easy to see that this vertex-splitting does not
change the sums of the weights of the matchings, as long as the
new edges we add are given weight 1.

Now we apply urban renewal in $n^2$ locations 
to get a graph of the form
$$
\pspicture(-3.5,-3.5)(3.5,3.5)
\pspolygon(-3,-2)(2,3)(3,2)(-2,-3)
\pspolygon(-3,0)(0,3)(3,0)(0,-3)
\pspolygon(-3,2)(-2,3)(3,-2)(2,-3)
\psline(-3.5,2)(-3,2)\psline(3,2)(3.5,2)
\psline(-3.5,0)(-3,0)\psline(3,0)(3.5,0)
\psline(-3.5,-2)(-3,-2)\psline(3,-2)(3.5,-2)
\psline(2,-3.5)(2,-3)\psline(2,3)(2,3.5)
\psline(0,-3.5)(0,-3)\psline(0,3)(0,3.5)
\psline(-2,-3.5)(-2,-3)\psline(-2,3)(-2,3.5)
\multips(-3,-2)(0,2){3}{\multips(0,0)(2,0){4}{\odisk(0,0){.1}}}
\multips(-2,-3)(0,2){4}{\multips(0,0)(2,0){3}{\odisk(0,0){.1}}}
\odisk(-3.5,2){.1}\odisk(-3.5,0){.1}\odisk(-3.5,-2){.1}
\odisk( 3.5,2){.1}\odisk( 3.5,0){.1}\odisk( 3.5,-2){.1}
\odisk(2,-3.5){.1}\odisk(0,-3.5){.1}\odisk(-2,-3.5){.1}
\odisk( 2,3.5){.1}\odisk( 0,3.5){.1}\odisk( -2,3.5){.1}
\endpspicture
$$
with new weights $w'$, $x'$, $y'$, $z'$ 
replacing the old weights $w$, $x$, $y$, $z$.
(If one were attempting a literal description of the algorithm,
one would want the weight-variables to have subscripts 
$i,j$ ranging from $1$ to $n$.)
Here is the coup de grace: The pendant edges that we see must
belong to \emph{every} perfect matching of the graph, 
so we can delete them from the graph, 
obtaining an Aztec graph of order $n-1$:
$$
\pspicture(-2,-2)(2,2)
\pspolygon(-2,-1)(1,2)(2,1)(-1,-2)
\pspolygon(-2,1)(-1,2)(2,-1)(1,-2)
\multips(-2,-1)(0,2){2}{\multips(0,0)(2,0){3}{\odisk(0,0){.1}}}
\multips(-1,-2)(0,2){3}{\multips(0,0)(2,0){2}{\odisk(0,0){.1}}}
\endpspicture
$$

The theorem now follows.

One can prove the ``fortress" theorem (see \cite{Ci2})
by similar means.  
Here, the initial graph is
$$
\pspicture(-2.5,-2.5)(2.5,2.5)
\multips(-1,-1.5)(0,1.5){3}{\multips(0,0)(1.5,0){2}{\psline(0,0)(.5,0)}}
\multips(-1.5,-1)(1.5,0){3}{\multips(0,0)(0,1.5){2}{\psline(0,0)(0,.5)}}
\psline(-2.5,-1.5)(-2,-1.5)\psline(-1.5,-2.5)(-1.5,-2)
\psline( 2.5,-1.5)( 2,-1.5)\psline( 1.5,-2.5)( 1.5,-2)
\psline(-2.5, 1.5)(-2, 1.5)\psline(-1.5, 2.5)(-1.5, 2)
\psline( 2.5, 1.5)( 2, 1.5)\psline( 1.5, 2.5)( 1.5, 2)
\multips(-1.5,-1.5)(0,1.5){3}{\multips(0,0)(1.5,0){3}{
    \pspolygon(0,.5)(.5,0)(0,-.5)(-.5,0)
    \odisk(0,.5){.1}\odisk(.5,0){.1}\odisk(0,-.5){.1}\odisk(-.5,0){.1}}}
\odisk(-2.5,-1.5){.1}\odisk(-1.5,-2.5){.1}
\odisk( 2.5,-1.5){.1}\odisk( 1.5,-2.5){.1}
\odisk(-2.5, 1.5){.1}\odisk(-1.5, 2.5){.1}
\odisk( 2.5, 1.5){.1}\odisk( 1.5, 2.5){.1}
\endpspicture
$$
One's first impulse is to remove the pendant edges, 
but it is better still to apply urban renewal 
to every \emph{other} ``city" in this graph.  
Then one gets an Aztec diamond graph in which 
roughly half of the edges have weight $1$ 
and the rest have weight $1/2$,
with the two sorts of cities alternating
in checkerboard fashion. 
The same techniques that were used above can also be applied here 
to allow one to conclude that the number of matchings of such a
graph is a certain power of 5 (or twice a certain power of $5$).  
The ``5" comes from the cell-factor $(1)(1)+(1/2)(1/2)$ 
(when multiplied by $4$), just as (in the unweighted Aztec diamond theorem) 
the ``2" comes from the cell-factor $(1)(1)+(1)(1)$.
In short, the proof of Yang's theorem becomes a near-triviality
(just as it is for Ciucu's method \cite{Ci2}).

\section{Computing edge-inclusion probabilities: the proof}

In this section I will prove that the iterative scheme for
computing these probabilities that was described earlier does in
fact work (at least when all edge-weights are non-vanishing).  
To start things off we will need a lemma:
This lemma is a generalization of a proposition
originally conjectured by Alexandru Ionescu
in the context of Aztec diamond graphs,
and was proved by Propp in 1993 (private e-mail communication).

Let $G$ be a bipartite planar graph
and let $A$,$B$,$C$,$D$ be four vertices of $G$
that form a 4-cycle $ABDC$ bounding a face of $G$:
$$
\pspicture(-2.5,-2.5)(2.5,2.5)
\pspolygon(-2,0)(0,2)(2,0)(0,-2)
\psline(2,0)(2.4,0)\psline(-2,0)(-2.4,0)
\psline(0,2)(0,2.4)\psline(0,-2)(0,-2.4)
\psline(.2,2.4)(0,2)(-.2,2.4)
\psline(.2,-2.4)(0,-2)(-.2,-2.4)
\psline(2.4,.2)(2,0)(2.4,-.2)
\psline(-2.4,.2)(-2,0)(-2.4,-.2)
\odisk(0,-2){.1} \odisk(0,2){.1}
\odisk(-2,0){.1} \odisk(2,0){.1}
\rput[r](-.2,2){$A$}
\rput[b](-2,.2){$B$}
\rput[b](2,.2){$C$}
\rput[r](-.2,-2){$D$}
\endpspicture
$$
Assume the edges of this graph are weighted in some fashion.
Define $G_{AB}$ as the weighted graph obtained 
by deleting vertices $A$ and $B$ 
and all incident edges;
all other edges have the same weight as in $G$. 
Define $G_{AC}$, $G_{BD}$, $G_{CD}$, 
and $G_{ABCD}$ similarly.  
For any edge-weighted graph $H$, let $M(H)$ denote 
the sum of the weights of the matchings of $H$.
Then the identity we will prove is
$$M(G) M(G_{ABCD}) = M(G_{AB}) M(G_{CD}) + M(G_{AC}) M(G_{BD}).$$

Note that if we superimpose a matching of $G$ 
and a matching of $G_{ABCD}$, 
we get a multiset of edges of $G$ 
in which vertices $A$,$B$,$C$,$D$ get degree 1 and
every other vertex gets degree 2.  
The same thing happens if we superimpose
a matching of $G_{AB}$ and a matching of $G_{CD}$, 
or a matching of $G_{AC}$ and a matching of $G_{BD}$.  
Call an edge-multiset of this sort a ``near 2-factor" of $G$.  
I will prove that equality holds above by showing that 
each near 2-factor $F$ makes the same contribution 
to the left side as it does to the right side.

Pick a specific near 2-factor $F$ of $G$; 
it consists of $k$ closed loops plus some paths 
that connect up the four vertices 
$A$,$B$,$C$,$D$ in pairs.  
It is impossible that the vertices 
are paired as $\{A,D\}$ and $\{B,C\}$, 
since the graph is planar.
Hence the vertices pair as 
either $\{A,B\}$ and $\{C,D\}$ 
or $\{A,C\}$ and $\{B,D\}$.  
Without loss of generality, 
let us focus on the former case.  
Since $G$ is bipartite, the path from $A$ to $B$ 
will contain an odd number of edges, 
as will the path from $C$ to $D$.
It is easy to verify
the following three claims:
\begin{enumerate}
\item There are $2^k$ ways to split up $F$ 
as a matching of $G$ together with 
a matching of $G_{ABCD}$.
\item There are $2^k$ ways to split up $F$ 
as a matching of $G_{AB}$ together with 
a matching of $G_{CD}$.
\item There are \emph{no} ways to split up $F$ 
as a matching of $G_{AC}$ together with 
a matching of $G_{BD}$.
\end{enumerate}
Moreover, one can check that each near 2-factor 
that pairs the vertices as $\{A,B\}$ and $\{C,D\}$
contributes equal weight to both 
sides of the equation.  
Clearly the same is true for the near 2-factors
that pair the vertices the other way
($A$ with $C$ and $B$ with $D$).
This completes the proof. 

(Note: Another way to prove this lemma
is to use the face that the number of perfect matchings
of a bipartite plane graph with $2n$ vertices
can be written as the determinant of an $n$-by-$n$
Kasteleyn-Percus matrix \cite{Ka} \cite{Pe}.
Then the lemma is seen to be a restatement
of Dodgson condensation (\cite{D}; see also \cite{RR}).

It will be convenient for us to restate this lemma
in a somewhat different form.
Assume that the graph has
at least one matching of positive weight.
Let $w$, $x$, $y$, and $z$ denote the respective weights
of the edges $AB$, $AC$, $BD$, and $CD$,
and let $p$, $q$, $r$, and $s$ denote the probabilities
of these respective edges being included
in a random matching of $G$
(where as usual the probability of an individual matching
is proportional to its weight):
$$
\pspicture(-2.5,-2.5)(2.5,2.5)
\pspolygon(-2,0)(0,2)(2,0)(0,-2)
\psline(2,0)(2.4,0)\psline(-2,0)(-2.4,0)
\psline(0,2)(0,2.4)\psline(0,-2)(0,-2.4)
\psline(.2,2.4)(0,2)(-.2,2.4)
\psline(.2,-2.4)(0,-2)(-.2,-2.4)
\psline(2.4,.2)(2,0)(2.4,-.2)
\psline(-2.4,.2)(-2,0)(-2.4,-.2)
\rput[br](-1.05, 1.05){$w\;\;[p]$} \rput[bl](1.05, 1.05){$x\;\;[q]$}
\rput[tr](-1.05,-1.05){$y\;\;[r]$} \rput[tl](1.05,-1.05){$z\;\;[s]$}
\odisk(0,-2){.1} \odisk(0,2){.1}
\odisk(-2,0){.1} \odisk(2,0){.1}
\rput[r](-.2,2){$A$}
\rput[b](-2,.2){$B$}
\rput[b](2,.2){$C$}
\rput[r](-.2,-2){$D$}
\endpspicture
$$
The probability $P^*$ that a random matching of $G$ 
has an alternating cycle at face $ABDC$
(i.e., either includes edges $AB$ and $CD$
or includes edges $AC$ and $BD$)
is equal to $M(G_{ABCD})(wz+xy)/M(G)$.
Since we have
$p=wM(G_{AB})/M(G)$,
$q=xM(G_{AC})/M(G)$,
$r=yM(G_{BD})/M(G)$, and
$s=zM(G_{CD})/M(G)$,
the preceding lemma tells us that,
in the event that $w$,$x$,$y$,$z$ are all non-zero,
the probability $P^*$ 
is equal to
$(psxy+qrwz)(\frac{1}{xy}+\frac{1}{wz})$.
It follows from this that
the probability that a random matching of $G$ contains
both edge $AB$ and edge $CD$ is $(psxy+qrwz)/xy$, 
while the corresponding probability for edges $AC$ and $BD$ 
is $(psxy+qrwz)/wz$.
(To see that this follows from the formula for $P^*$,
note that their sum is indeed $(wz+xy)(ps/wz+qr/xy)=P^*$,
and that their ratio is $wz/xy$, 
which is the ratio of the weights of any two matchings 
that differ only in that one contains $AB$ and $CD$ 
while the other contains $AC$ and $BD$.)

Now let us go back and give urban renewal a fresh look.  Here are
the graphs $G$ and $G'$, with edges marked with both their weight and
their probability (in the format ``weight [probability]").
Here is $G$:
$$
\pspicture(-2.5,-2.5)(2.5,2.5)
\pspolygon(-1,0)(0,1)(1,0)(0,-1)
\psline(1,0)(2.4,0)\psline(-1,0)(-2.4,0)
\psline(0,1)(0,2.4)\psline(0,-1)(0,-2.4)
\psline(.2,2.4)(0,2)(-.2,2.4)
\psline(.2,-2.4)(0,-2)(-.2,-2.4)
\psline(2.4,.2)(2,0)(2.4,-.2)
\psline(-2.4,.2)(-2,0)(-2.4,-.2)
\rput[br](-.55, .55){$w\;\;[p]$} \rput[bl](.55, .55){$x\;\;[q]$}
\rput[tr](-.55,-.55){$y\;\;[r]$} \rput[tl](.55,-.55){$z\;\;[s]$}
\odisk(0,-2){.1} \odisk(0,-1){.1} \odisk(0, 1){.1} \odisk(0, 2){.1}
\odisk(-2,0){.1} \odisk(-1,0){.1} \odisk( 1,0){.1} \odisk( 2,0){.1}
\rput[r](-.2,2){$A$}
\rput[b](-2,.2){$B$}
\rput[b](2,.2){$C$}
\rput[r](-.2,-2){$D$}
\endpspicture
$$
(The unmarked edges have weight $1$, and their probabilities are
$1-p-q$, $1-q-s$, et cetera.)

Here is $G'$:
$$
\pspicture(-2.5,-2.5)(2.5,2.5)
\pspolygon(-2,0)(0,2)(2,0)(0,-2)
\psline(2,0)(2.4,0)\psline(-2,0)(-2.4,0)
\psline(0,2)(0,2.4)\psline(0,-2)(0,-2.4)
\psline(.2,2.4)(0,2)(-.2,2.4)
\psline(.2,-2.4)(0,-2)(-.2,-2.4)
\psline(2.4,.2)(2,0)(2.4,-.2)
\psline(-2.4,.2)(-2,0)(-2.4,-.2)
\rput[br](-1.05, 1.05){$w'\;\;[p']$} \rput[bl](1.05, 1.05){$x'\;\;[q']$}
\rput[tr](-1.05,-1.05){$y'\;\;[r']$} \rput[tl](1.05,-1.05){$z'\;\;[s']$}
\odisk(0,-2){.1} \odisk(0,2){.1}
\odisk(-2,0){.1} \odisk(2,0){.1}
\rput[r](-.2,2){$A$}
\rput[b](-2,.2){$B$}
\rput[b](2,.2){$C$}
\rput[r](-.2,-2){$D$}
\endpspicture
$$
Recall that we are in the situation where we know $p'$, $q'$, $r'$, and $s'$,
and are trying to compute $p$, $q$, $r$, and $s$.  Recall also that
\begin{align*}
w' &= z/(wz+xy) & x' &= y/(wz+xy) \\
y' &= x/(wz+xy) & z' &= w/(wz+xy).
\end{align*}

Here is our plan: First, we will work out the probabilities of all the
local patterns in $G'$.  Then we will use the urban renewal correspondence
to deduce the probabilities of the local patterns in $G$.  Finally,
we will deduce the probabilities of the individual edges in $G$.

I will write $p'$,$q'$,$r'$,$s'$,$w'$,$x'$,$y'$,$z'$ as
$P$,$Q$,$R$,$S$,$W$,$X$,$Y$,$Z$, to save eye-strain.  
I will also write $\Delta'=XYPS+WZQR$.

Let $S$ denote the set of vertices in $\{A,B,C,D\}$ that are matched to
a vertex inside the patch.  Then here are the respective probabilities
of the local patterns in $G'$:
\begin{align*}
S &= \{A,B,C,D\}: (WZ+XY)\Delta'/WXYZ \\
S &= \{A,B\}:     P-\Delta'/XY \\
S &= \{A,C\}:     Q-\Delta'/WZ \\
S &= \{B,D\}:     R-\Delta'/WZ \\
S &= \{C,D\}:     S-\Delta'/XY \\
S &= \{\}:        1-P-Q-R-S+(WZ+XY)\Delta'/WXYZ
\end{align*}
(The second formula is obtained by recalling that the probability of
edge $AB$ being present in a random matching, which can happen in two
ways according to whether or not edge $CD$ is present, must be $P$;
the last formula is obtained by recalling that the probabilities of
the six local configurations must sum to 1.)

It follows from the last of these, and the urban renewal correspondence,
that the probability that a random matching of $G$ contains two edges
from the 4-cycle $ABCD$ is
$$1-P-Q-R-S+(WZ+XY)(PS/WZ+QR/XY).$$
The probability that a random matching of $G$ 
contains edges $AB$ and $CD$
must be equal to this quantity times
$$wz/(wz+xy) = WZ/(WZ+XY),$$
which gives
$$(WZ/(WZ+XY))(1-P-Q-R-S) + PS + QR(WZ/XY).$$
The probability that a random matching of $G$ contains edge $AB$ 
but \emph{not} edge $CD$ 
(by another application of the urban renewal correspondence)
equals the probability that a random matching of $G'$ 
contains edge $CD$ but not edge $AB$, 
which is $$S - PS - QR(WZ/XY).$$
Adding, we find that the probability that a random matching of $G$ 
contains edge $AB$ is
$$S + \frac{wz}{wz+xy}(1-P-Q-R-S).$$
Thus we conclude that
$$
p = S + \frac{wz}{wz+xy} (1 - P - Q - R - S) ,
$$
where $wz/(wz+xy)$ is the creation bias and
$1 - P - Q - R - S$ is the deficit.
The formulas for $q$, $r$, and $s$ follow by symmetry.

\section{Generating a random matching: the proof}

Now I will show that (generalized) domino-shuffling works.

First, let us note that the algorithm for generating
random matchings that was presented earlier
really \emph{is} a generalization of the domino-shuffling procedure
described in \cite{EKLP}.  
The operation of removing from a matching those matched edges 
that belong to the same cell as another matched edge 
corresponds to the process of removing odd blocks 
to obtain an odd-deficient tiling;
the operation of swapping the remaining edges 
to the opposite side of the cell 
corresponds to the process of shuffling dominos 
(resulting in the creation of an odd-deficient tiling 
of a larger Aztec diamond);
and the operation of introducing new edges 
corresponds to the creation
of new $2 \times 2$ blocks.

Next, let us recall that a perfect matching
of the Aztec diamond graph of order $n$
determines an $n \times n$ alternating-sign matrix,
as was first described in \cite{EKLP}.
Specifically, for each of the $n^2$ cells of the graph, 
write down the number of edges from that cell 
that participate in the matching and subtract 1.  
For instance:
$$
\pspicture(-4,-4)(4,4)
\psline(-4,3)(-3,4)\psline(-2,3)(-1,4)\psline(1,4)(2,3)\psline(3,4)(4,3)
\psline(-1,2)(0,3)\psline(-3,2)(-2,1)\psline(1,2)(2,1)\psline(3,2)(4,1)
\psline(-4,1)(-3,0)\psline(0,1)(1,0)\psline(-2,-1)(-1,0)\psline(2,-1)(3,0)
\psline(-4,-1)(-3,-2)\psline(0,-1)(1,-2)\psline(3,-2)(4,-1)\psline(-2,-3)(-1,-2)
\psline(-4,-3)(-3,-4)\psline(-1,-4)(0,-3)\psline(1,-4)(2,-3)\psline(3,-4)(4,-3)
\multips(-4,-3)(0,2){4}{\multips(0,0)(2,0){5}{\odisk(0,0){.1}}}
\multips(-3,-4)(0,2){5}{\multips(0,0)(2,0){4}{\odisk(0,0){.1}}}
\rput(-3, 3){$ 0$}\rput(-1, 3){$+1$}\rput(1, 3){$ 0$}\rput(3, 3){$ 0$}
\rput(-3, 1){$+1$}\rput(-1, 1){$-1$}\rput(1, 1){$+1$}\rput(3, 1){$ 0$}
\rput(-3,-1){$ 0$}\rput(-1,-1){$ 0$}\rput(1,-1){$ 0$}\rput(3,-1){$+1$}
\rput(-3,-3){$ 0$}\rput(-1,-3){$+1$}\rput(1,-3){$ 0$}\rput(3,-3){$ 0$}
\endpspicture
$$
Label the weights of the edges of the cell in row $i$ and column $j$ 
with $w_{ij}$, $x_{ij}$, $y_{ij}$, and $z_{ij}$, thus:
$$
\pspicture(-1,-1)(1,1)
\pspolygon(-1,0)(0,1)(1,0)(0,-1)
\rput[br](-.55, .55){$w$} \rput[bl](.55, .55){$x$}
\rput[tr](-.55,-.55){$y$} \rput[tl](.55,-.55){$z$}
\odisk(0,-1){.1} \odisk(0, 1){.1}
\odisk(-1,0){.1} \odisk( 1,0){.1}
\endpspicture
$$
Define also a ``cell-weight"
$$D_{ij} = w_{ij} z_{ij} + x_{ij} y_{ij}$$
(these weights are the same as the cell-factors
considered earlier in this article).
Call two edges in a matching ``parallel" 
if they belong to the same cell.

The bias in the creation rates for generalized shuffling 
was chosen so that the ratio between
the probabilities of choosing 
one pair of parallel edges versus another
($wz$ versus $xy$)
is equal to the ratio of the weights
of the resulting matchings.
Consequently, if two perfect matchings 
of the Aztec graph of order $n$ 
are associated with the same $n \times n$ alternating-sign matrix, 
then their \emph{relative} probabilities are correct.  
So all we need to do is verify that the aggregate probability 
of these matchings (for fixed alternating-sign matrix $A$) 
is what it should be.  

To do this, it is convenient to introduce
certain partial matchings of Aztec diamond graphs
(analogous to the odd-deficient and even-deficient
domino-tilings considered in \cite{EKLP}).
For the Aztec diamond of order $n$,
we consider partial matchings obtained from perfect matchings
by removal of all parallel edges;
that is, pairs of edges that
belong to the same cell.
For the Aztec diamond of order $n-1$,
we consider partial matchings obtained from perfect matchings
by removal of all pairs of edges
that belong to the same ``co-cell''
(defined as the space between four adjoining cells).
The ``destruction'' phase of generalized shuffling
converts a perfect matching of the graph of order $n-1$
into a partial matching $M$ of the graph of order $n-1$,
the ``sliding'' phase converts that partial matching
into a partial matching $M'$ of the graph of order $n$,
and the ``creation'' phase converts that partial matching
into a perfect matching.
Note that $M$ determines $M'$ uniquely, and vice versa.

Define the aggregate cell-weights $D_+$, $D_0$, and $D_-$
to be the products of those cell-weights $D_{i,j}$
associated with locations in the Aztec diamond
at which $+1$, $0$, and $-1$ appear, respectively
(so that $D_+ D_0 D_-$ is just the product of all the cell-weights).
The following three claims can be readily checked:
the sum of the weights of the perfect matchings
of the Aztec diamond of order $n$ that extend the partial matching $M'$ 
is equal to $D_+$ times the weight of $M'$ itself
(defined as the product of the weights of its constituent edges);
the weight of $M'$ is equal to $D_0$ times the weight of $M$
(defined in terms of the weighting of the graph of order $n-1$);
and the weight of $M$ is equal to $D_-$ 
times the sum of the weights of the perfect matchings
of the Aztec diamond of order $n-1$ that extend the partial matching $M$. 
(This relationship between cell-weights
and the entries of alternating-sign matrices
was originally made clear in the work of Ciucu \cite{Ci1}.)

The upshot is, the sum of the weights of the perfect matchings
of the Aztec diamond graph of order $n$ that are associated
with the alternating-sign matrix $A$ equals
$$
D_{1,1} D_{1,2} \ldots D_{n,n}
$$
times the sum of the weights of all the perfect matchings
of the Aztec diamond graph of order $n-1$
that are capable of giving rise to it under shuffling.
This factor is independent of $A$.
It follows that if one takes a random matching of the smaller graph
(in accordance with the edge weights given by urban renewal) and
applies destruction, shuffling, and creation, then one will get a random
matching of the larger graph.

\section{Diabolo-tilings of fortresses}

Consider an $n$-by-$n$ array of unit squares,
each of which has been cut by both of its diagonals, 
so that there are $4n^2$ isosceles right triangles in all.
A region called a \emph{fortress of order n}
is obtained by removing some of the triangles
that are at the boundary of the array.
More specifically, one colors the triangles alternately black and white,
and removes all the black (resp.\ white) triangles
that have an edge on the top or bottom (resp.\ left or right) 
boundaries of the array. 
Here, for instance, is a fortress of order 4:
$$
\pspicture(-4,-4)(4,4)
\psset{xunit=.8cm,yunit=.8cm}
\psline(-4,4)(-2,4)(-1,3)(0,4)(2,4)(3,3)
\psline(4,-4)(2,-4)(1,-3)(0,-4)(-2,-4)(-3,-3)
\pspolygon(-4,2)(-2,4)(4,-2)(2,-4)
\pspolygon(4,2)(2,4)(-4,-2)(-2,-4)
\pspolygon(4,0)(0,4)(-4,0)(0,-4)
\multips(-4,0)(2,-2){2}{\pspolygon[fillstyle=solid,fillcolor=black](0,0)(1,-1)(0,-2)}
\multips(-4,4)(2,-2){4}{\pspolygon[fillstyle=solid,fillcolor=black](0,0)(1,-1)(0,-2)}
\multips(0,4)(2,-2){2}{\pspolygon[fillstyle=solid,fillcolor=black](0,0)(1,-1)(0,-2)}
\multips(-3,-1)(2,-2){2}{\pspolygon[fillstyle=solid,fillcolor=black](0,0)(1,1)(1,-1)}
\multips(-3,3)(2,-2){4}{\pspolygon[fillstyle=solid,fillcolor=black](0,0)(1,1)(1,-1)}
\multips(1,3)(2,-2){2}{\pspolygon[fillstyle=solid,fillcolor=black](0,0)(1,1)(1,-1)}
\multips(-4,2)(2,2){1}{\pspolygon[fillstyle=solid,fillcolor=black](0,0)(2,0)(1,-1)}
\multips(-4,-2)(2,2){3}{\pspolygon[fillstyle=solid,fillcolor=black](0,0)(2,0)(1,-1)}
\multips(0,-2)(2,2){2}{\pspolygon[fillstyle=solid,fillcolor=black](0,0)(2,0)(1,-1)}
\multips(-4,0)(2,2){2}{\pspolygon[fillstyle=solid,fillcolor=black](0,0)(2,0)(1,1)}
\multips(-2,-2)(2,2){3}{\pspolygon[fillstyle=solid,fillcolor=black](0,0)(2,0)(1,1)}
\multips(2,-2)(2,2){1}{\pspolygon[fillstyle=solid,fillcolor=black](0,0)(2,0)(1,1)}
\endpspicture
$$
(When $n$ is even, the two ways of coloring the triangles
lead to fortresses that are mirror images of one another;
however, when $n$ is odd, the two sorts of fortresses
one obtains by removing the colored triangles in the
specified fashion are genuinely different.)
The graph dual to a fortress 
(with vertices corresponding to the isosceles right triangles)
is composed of 4-cycles and 8-cycles.
(If one removes the pendant edges
from the final figure in section 5
one obtains the graph dual to a fortress of order 3.)

We call the small isosceles right triangles \emph{monobolos}.
Given our tiling of the fortress by monobolos,
we can (in many different ways) form a new tiling
with half as many tiles
by amalgamating the monobolos in pairs,
where the two monobolos that are paired together
are required to be adjacent.
These new tiles are \emph{diabolos},
and come in two different shapes:
squares and isosceles triangles.
(Note: in the literature on tilings, such as \cite{G},
a third kind of diabolo is recognized,
namely a parallelogram with angles
of 45 and 135 degrees with side-lengths
in the ratio $1:\sqrt{2}$;
however, these diabolos do not arise
in the context being discussed here.)

A tiling obtained from the original monobolo-tiling of the fortress
by amalgamating adjacent pairs of monobolos
will be called, for (comparative) brevity,
a diabolo-tiling of a fortress.
It should, however, be kept in mind that we are tacitly
limiting ourselves to those tilings that can be obtained
by the aforementioned process of amalgamation;
or, if one prefers to use the dual picture,
tilings that correspond to perfect matchings
of the ``squares and octagons'' graph that is dual 
to the original tiling by monobolos.
For instance, the first of the following two tilings
of a fortress of order 3
by means of tiles that are diabolos
is an example of what we mean by a
``diabolo-tiling of a fortress'',
but the second is not,
because the four center diabolos
are not amalgams of the required kind:
$$
\pspicture(-3,-3)(3,3)
\psset{xunit=.8cm,yunit=.8cm}
\psline(-1,3)(1,3)
\psline(1,1)(3,1)
\psline(-3,-1)(-1,-1)
\psline(-1,-3)(1,-3)
\psline(-3,1)(-3,-1)
\psline(-1,-1)(-1,-3)
\psline(1,3)(1,1)
\psline(3,1)(3,-1)
\psline(-3,1)(-1,3)
\psline(-2,0)(0,2)
\psline(0,-2)(2,0)
\psline(1,-3)(3,-1)
\psline(-3,-1)(-1,-3)
\psline(-3,1)(1,-3)
\psline(-2,2)(-1,1)
\psline(1,-1)(2,-2)
\psline(-1,3)(3,-1)
\psline(1,3)(3,1)
\psline(-1,-1)(1,1)
\psline(-1,1)(1,-1)
\endpspicture
$$
$$
\pspicture(-3,-3)(3,3)
\psset{xunit=.8cm,yunit=.8cm}
\psline(-1,3)(1,3)
\psline(1,1)(3,1)
\psline(-3,-1)(-1,-1)
\psline(-1,-3)(1,-3)
\psline(-3,1)(-3,-1)
\psline(-1,-1)(-1,-3)
\psline(1,3)(1,1)
\psline(3,1)(3,-1)
\psline(-3,1)(-1,3)
\psline(-2,0)(0,2)
\psline(0,-2)(2,0)
\psline(1,-3)(3,-1)
\psline(-3,-1)(-1,-3)
\psline(-3,1)(1,-3)
\psline(-2,2)(-1,1)
\psline(1,-1)(2,-2)
\psline(-1,3)(3,-1)
\psline(1,3)(3,1)
\psline(-2,0)(2,0)
\psline(0,2)(0,-2)
\endpspicture
$$

To count the diabolo-tilings of a fortress
(understood in the above sense),
one replaces the underlying monobolo-tiling
by its dual graph and adds pendant edges
along the boundary,
obtaining a graph like the one shown
at the end of section 5.
One can then apply the urban renewal lemma
to turn this into an Aztec diamond graph
with edges of weight 1 arranged in cells
that alternate, checkerboard-style,
with cells whose edges have weight $\frac12$.
In this fashion one can derive
Yang's formula.

Now, however, we set our sights higher.
In section 6, 
we gave a scheme for calculating
the probability that a given edge of a
weighted Aztec diamond graph
is included in a random matching,
where the probability of a matching
is proportional to its weight.
This in fact lets us write down an explicit power series
whose coefficients are precisely
those inclusion-probabilities.
We will show how this goes in a
slightly more general setting,
where the edges of weight $\frac12$
are replaced by edges of weight $t$. 

Here we will find it convenient to use 
a rotated and shifted version of the Aztec diamond graph of order $n$, 
whose edges are all either horizontal or vertical (rather than diagonal).
The vertices of this graph are the integer points $(i,j)$
where $|i+\frac12|+|j+\frac12| \leq n$;
the edges of the graph join vertices at distance 1.
The horizontal edge joining $(i-1,j)$ and $(i,j)$
is called ``north-going'' if $i+j+n$ is odd
and ``south-going'' otherwise,
and the vertical edge joining $(i,j-1)$ and $(i,j)$
is called ``east-going'' if $i+j+n$ is odd
and ``west-going'' otherwise
(this nomenclature is taken from \cite{EKLP}).
The kind of weighting we are considering can be described as follows:
\begin{enumerate}
\item
Each horizontal edge has weight 1 or weight $t$ according to the parity
of its $x$-coordinate.  (That is, if two such edges are related by
a unit vertical displacement, they have same weight, but if they 
are related by a unit horizontal displacement, their weights differ.)
\item
Each vertical edge has weight 1 or weight $t$ according to the parity
of its $y$-coordinate.
\item
If $n$ is 1 or 2 (mod 4), the four extreme-most edges (the northernmost
and southernmost horizontal edges and the easternmost and westernmost
vertical edges) have weight $t$; if $n$ is 0 or 3 (mod 4), those four edges 
have weight 1.
\end{enumerate}
(These four situations are considered together because
each is linked to the next by the shuffling operation.
Technically, this scheme only handles fortresses of odd order,
and indeed, only handles half of those,
but the other case are quite similar.)
We define a generating function $P(x,y,z)$ in which
the coefficient of $x^i y^j z^n$
(with $k \geq 0$, $i+j+k$ odd, and $|i|+|j|<n$)
is the probability that a random matching 
of the Aztec diamond graph of order $n$, 
chosen in accordance with the probability distribution 
determined by the weight $t$, contains 
a horizontal bond joining vertices $(i-1,j)$ and $(i,j)$,
where $i+j+n$ is odd. 
(We are limiting ourselves to northgoing horizontal bonds
but this limitation is unimportant,
since by symmetry the other three sorts of bonds
have the same behavior.)
Note that the coefficients in this generating function
(viewed as a power series in the variables $x$, $y$, and $z$)
are polynomials in $t$,
and that for each $n$,
only finitely many pairs $i,j$
contribute a non-zero term to $P(x,y,z)$.

The presence of $t$ in the pattern of edge-weighting
makes it natural to divide the north-going edges
into four classes, according to how they ``sit''
relative to the weighting
(cells are weighted in two patterns when $n$ is odd
and in two other patterns when $n$ is even).
Correspondingly, we write $P(x,y,z)$ as a sum
of four generating functions,
each corresponding to one of the four classes of north-going edges.
We call these generating functions $E$, $F$, $G$, and $H$.
We also introduce generating functions whose coefficients
are the net creation rates of the four kinds 
of differently-weighted cells,
which we denote by $e$, $f$, $g$, and $h$.
$E$ and $e$ are associated with cells 
in which all edges have weight $t$;
$F$ and $f$ are associated with cells 
in which all edges have weight 1;
$G$ and $g$ are associated with cells 
in which horizontal edges have weight $t$ and vertical edges have weight 1;
and $H$ and $h$ are associated with cells 
in which horizontal edges have weight 1 and vertical edges have weight $t$.
Generalized shuffling implies the following algebraic relationships
among the four generating functions:
\begin{align*}
e &= 1 + \frac{t^2}{1+t^2} z ( (x+x^{-1}) g + (y+y^{-1}) h ) - z^2 f \\
f &=       \frac{1}{1+t^2} z ( (x+x^{-1}) h + (y+y^{-1}) g ) - z^2 e \\
g &=          \frac12      z ( (x+x^{-1}) f + (y+y^{-1}) e ) - z^2 h \\
h &=          \frac12      z ( (x+x^{-1}) e + (y+y^{-1}) f ) - z^2 g
\end{align*}
These can be solved simultaneously, giving us four three-variable
generating functions $e(x,y,z)$, $f(x,y,z)$, $g(x,y,z)$, $h(x,y,z)$.
We can then solve
\begin{align*}
E &= Hyz + \frac{1}{2}       ez &
F &= Gyz + \frac{1}{2}       fz \\
G &= Eyz + \frac{t^2}{1+t^2} gz &
H &= Fyz + \frac{1}{1+t^2}   hz
\end{align*}
and obtain (messy!) generating functions $E(x,y,z)$, $F(x,y,z)$,
$G(x,y,z)$, $H(x,y,z)$.  The sum $E+F+G+H$, multiplied by $z$,
is the desired generating function $P$.  It is too messy to write 
down here, but it does at least in principle give one a way to 
understand what is going on with random diabolo-tilings of large 
fortresses.

We can also use generalized shuffling to generate a random tiling
of a large fortress.  Recall that when this was done for Aztec
diamonds in the early 1990s, a new phenomenon was discovered:
the ``arctic circle'' effect \cite{JPS}.  Specifically, with 
probability going to 1 as the size of the Aztec diamond graph
increases, a random tiling divides itself up naturally into five 
regions: four outer, ``frozen'' regions in which all the dominos 
line up with their neighbors in a repeating pattern, and a central,
``temperate'' region in which the dominos are jumbled together
in a random-looking way.  One feature of the temperate zone,
rigorously proved in \cite{CEP}, is that the asymptotic local 
statistics, expressed as a function of normalized position,
are nowhere constant; that is, the probability of seeing a
particular local pattern of dominos in a random tiling changes
as one shifts one view over macroscopic distances (i.e.,
distances comparable to the size of the Aztec diamond graph).
Strikingly, the boundary of the temperate zone, rescaled,
tends in probability to a perfect circle.

\begin{figure*}
\begin{center}
\scalebox{.9428}{\includegraphics{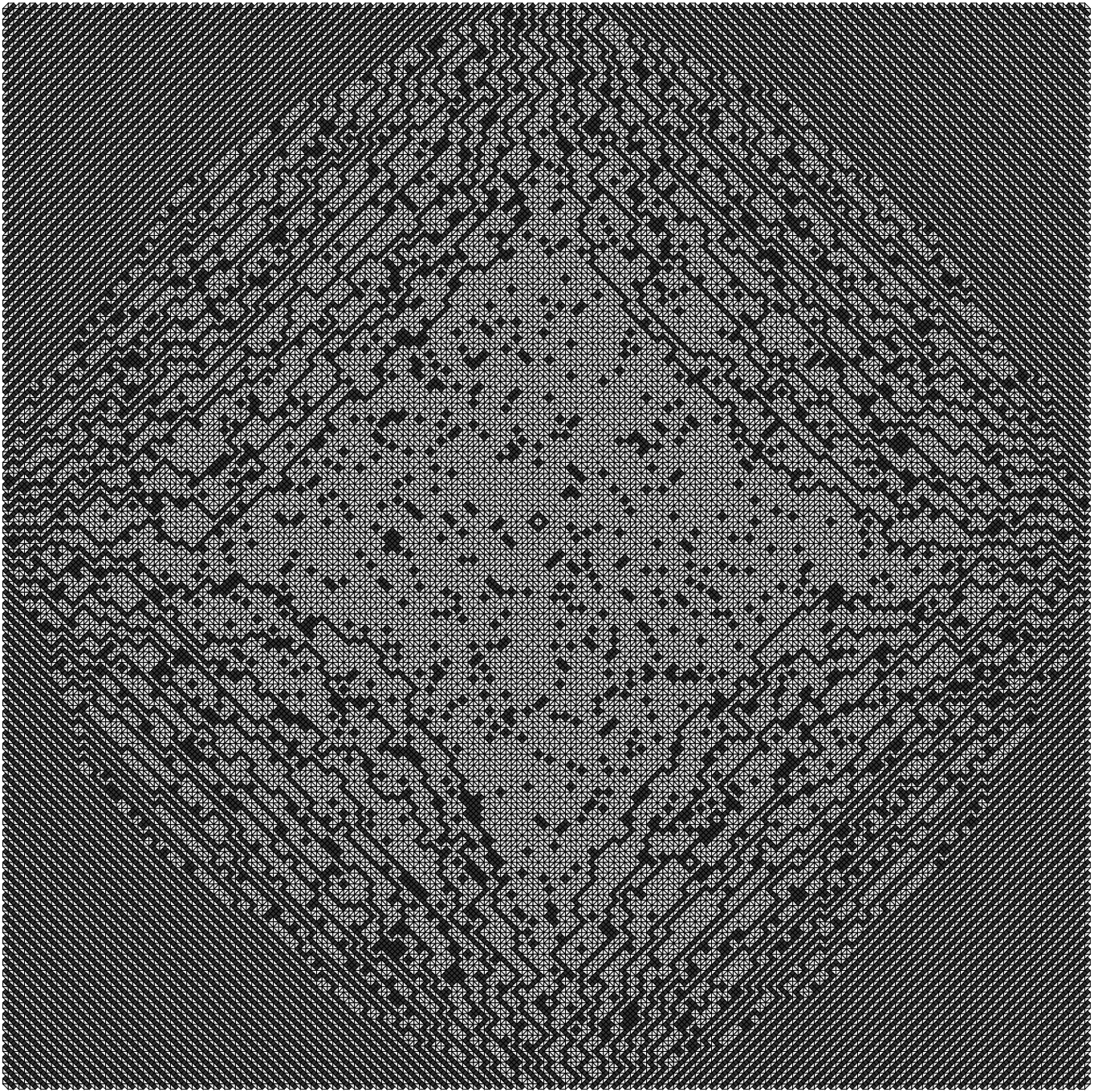}}
\caption{A fortress of order 200 randomly tiled with diabolos}
\label{figfort}
\end{center}
\end{figure*}

One might expect random diabolo-tilings of fortresses to exhibit
much the same sort of phenomena, albeit with the circular
temperate zone probably replaced by a temperate zone of
some other shape.  However, when the generalized shuffling
algorithm was used to generate random tilings of large
fortresses, a startling new phenomenon appeared: within the
very inmost part of the fortress, local statistics do \emph{not}
appear to undergo variation.  That is to say, within this
region, dubbed the ``tropical astroid'', the local statistics
of a random tiling appear to be shielded from the boundary,
so that the (normalized) position of a tile relative to the
boundary does not matter (as long as it stays fairly close
to the center of the region).  Figure \ref{figfort} shows a
random diabolo-tiling of the fortress of order 200; 
the square diabolos have been shaded, to highlight
the shape of the curve that separates the tropical zone
from the temperate zone.

Suppose we put coordinates on the fortress so that the corners are
at $(0,2)$, $(0,-2)$, $(2,0)$, and $(-2,0)$.  Then it appears
that the boundary of the frozen region
is given by one real component of the curve
\begin{multline*}
  400x^8 + 400y^8 +3400x^2y^6 + 3400y^2x^6 + 8025x^4y^4 \\
+ 1000x^6 + 1000y^6 -17250x^4y^2 -17250x^2y^4 -1431x^4 \\
-1431y^4 +25812x^2y^2 -3402x^2 -3402y^2 +729 = 0
\end{multline*}
(the ``octic circle'')
and that the boundary of the tropical region
is given by the other real component.
Henry Cohn and Robin Pemantle, as of this writing,
are working on a proof of these two assertions,
by giving an asymptotic analysis of the coefficients
of the generating functions
obtained via generalized shuffling. 

\section{Last thoughts}

Here are some last thoughts 
concerning urban renewal, domino-shuffling et cetera.

\begin{description}
\item[1.] This article has glossed over an important point,
namely, anomalies that arise when one attempts to apply
the algorithm to Aztec diamond graphs 
in which some edges have been assigned weight 0.  
Even if the weighted Aztec diamond in question has matchings
of positive weight,
it is possible that somewhere in the reduction process
one will encounter cells with cell-weight 0;
this leads to blow-up problems when one attempts to
divide by the cell-weight.
In such a case, one should replace edge-weights equal to $0$
in the original graph with edge-weights equal to $\epsilon$; 
the result will be a rational function of $\epsilon$ 
whose limiting value as $\epsilon$ goes to zero is desired.  
For this purpose, every rational function of $\epsilon$
can be written as a polynomial or power series in $\epsilon$ 
and replaced by its leading term, 
so calculations are not as bad as one might think.
This will work as long as the original weighted Aztec diamond graph
has at least one matching of positive weight.

\item[2.] These algorithms arose partly in response to work of
Kuperberg (see \cite{Ku1}), 
which took a more algebraic perspective on enumeration, 
using the approach pioneered by Kasteleyn \cite{Ka}.
Kasteleyn's method requires making some arbitrary choices
that end up not affecting the final answers
to meaningful enumerative questions.
This new work arose out of an attempt to find the invariant
combinatorial core of the Kasteleyn method.
In particular, graph-rewriting is a combinatorial substitute
for algebraic operations like row-reduction 
applied to a Kasteleyn matrix.
However, this analogy was never worked out in any kind
of rigorous detail.
It would be helpful to see this explained.

\item[3.] Continuing the above remark:
One strength of the algebraic approach,
exploited by Kenyon \cite{Ke1} and others,
is that edge-inclusion probabilities,
and, more generally, ``pattern-occurrence probabilities'',
can be expressed in terms of determinants
of minors of the inverse of the Kasteleyn matrix.
Here again, the answer must be independent of
the arbitrary choices that were made
in forming the matrix.
So it is natural to hope that there will be
a similar canonical scheme for calculating
such pattern-occurrence probabilities as well.
This may be similar to the problem of 
finding an extension of Dodgson's condensation scheme \cite{D}
that permits one to efficiently compute the inverse of a matrix
and not just compute its determinant.

\item[4.] Generalized shuffling
bears a strong resemblance
to an algorithm described in recent work of Viennot \cite{V}.
I have not studied the matter deeply enough
to be able to identify the relationship,
but I strongly believe that the resemblance
is more than coincidental.

\item[5.] It would also be desirable to extend urban renewal to 
matchings of \emph{non}-bipartite planar graphs.  This is
equivalent to asking for routinized matrix-reduction schemes 
for Kasteleyn's Pfaffians.

\item[6.] As described in \cite{CEP},
Ionescu's recurrence gives rise to edge-inclusion probabilities
for uniformly-weighted perfect matchings of Aztec diamond graphs,
and it is observed that if two edges $e$, $e'$
are close to one another (relative to the overall dimensions
of the large Aztec diamond graph they both inhabit),
\emph{and} $e$ and $e'$ occupy identical positions
within their respective cells
(i.e., both are either the northwest, northeast,
southwest, or southeast edges
in their cells),
then the edge-inclusion probabilities for $e$ and $e'$
are very close numerically.
Indeed, this phenomenon is a rigorously-proved consequence
of the detailed analysis given in \cite{CEP}.
However, it would be good to have a conceptual explanation
for this continuity property.
In particular, it seems plausible that
the recurrence relation for edge-inclusion probabilities
has the property of smoothing out differences.
A clear formulation of such a smoothing-out property,
and a rigorous proof that it holds,
would be very desirable,
since it might lead to a proof
that this kind of continuity property holds
for lozenge tilings of hexagons.
(Numerical evidence supports this continuity conjecture,
but the method used in \cite{CLP}
does not permit one to draw conclusions of this nature.) 
On the other hand, it should be noted that
these issues are subtle;
for instance, net creation rates
(viewed as a function of position)
tend not to vary as smoothly
as edge-inclusion probabilities. 

\item[7.] As is explained in \cite{EKLP},
enumeration of perfect matchings of Aztec diamond graphs
corresponds to ``2-enumeration'' of alternating-sign matrices,
where the weight assigned to a particular ASM
is 2 to the power of the number of $-1$'s it contains.
One can more generally consider $x$-enumeration,
where 2 is replaced by a general quantity $x$;
ordinary enumeration corresponds to 1-enumeration,
and 3-enumeration of ASMs leads to
interesting exact formulas
(see e.g.\ \cite{Ku2}).
Could shuffling be extended to give a scheme for
sampling from the set of ASMs with uniform distribution
or more generally an $x$-weighted distribution?

\end{description}

\end{document}